\documentclass[12pt]{article}

\usepackage{amssymb}
\usepackage{amsmath}
\usepackage{amsbsy}
\usepackage{amscd}
\usepackage{amsfonts}
\usepackage{amsthm}
\usepackage{mathrsfs}
\usepackage{verbatim}
\usepackage[colorlinks]{hyperref}
\usepackage{fullpage}
\usepackage{mathdots}
\usepackage{graphicx,subfigure}
\usepackage[english]{babel}
\usepackage[utf8]{inputenc}
\usepackage{tikz}
\usepackage{stackrel}
\usepackage{authblk}
\usepackage{stmaryrd}
\usepackage{tikz-cd}

\usepackage{mathtext}

\usetikzlibrary{backgrounds,fit, matrix}
\usetikzlibrary{positioning}
\usetikzlibrary{calc,through,chains}
\usetikzlibrary{arrows,shapes,snakes,automata, petri}

\usepackage{booktabs}
\usepackage{adjustbox}

\newtheorem{Theorem}[equation]{Theorem}
\newtheorem{Corollary}[equation]{Corollary}
\newtheorem{Lemma}[equation]{Lemma}
\newtheorem{Proposition}[equation]{Proposition}

\theoremstyle{definition}
\newtheorem{Definition}[equation]{Definition}
\newtheorem{Example}[equation]{Example}

\newtheorem{Remark}[equation]{Remark}

\numberwithin{equation}{section}
\numberwithin{figure}{section}

\newcommand{\C}{{\mathbb C}}
\newcommand{\Z}{{\mathbb Z}}

\newcommand{\N}{{\mathbb N}}

\newcommand{\mt}[1]{\text{#1}}

\newcommand{\OO}{\textbf{\"O}}

\begin{document}

\title{Asymptotic Semigroups and Two-sided Weak Orders}

\author{Mahir Bilen Can}
\affil[1]{{mahirbilencan@gmail.com}}

\maketitle

\begin{abstract}

Various partial orders related to the structures of dual canonical monoids are investigated. 
It is shown that the nilpotent variety of a dual canonical monoid is equidimensional; its dimension is found.
It is shown in type A that certain intervals of the Putcha poset of a dual canonical monoid are 
isomorphic to the Renner monoids of matrices. 
The notion of a two-sided weak order on a normal reductive monoid is introduced.  
A criterion, in terms of type maps, for the covering relations in a two-sided weak order to have degree 2 is found. 
It is shown that, for the unique equivariant divisor of a dual canonical monoid (the asymptotic semigroup), 
the covering relations of the two-sided weak order are always of degree 1. 
These  computations provide new insights for the two-sided weak orders on Coxeter groups. 
In type A, some enumerative results for the covering relations are presented.

\vspace{.5cm}


\end{abstract}

\normalsize

\section{Introduction}

Let $M$ be a complex reductive monoid with zero. 
Let $M_{nil}$ denote the variety of nilpotent elements of $M$.
Let $G$ denote the unit group of $M$. Let $B$ be a Borel subgroup of $G$. 
The following finite decompositions of $M$ and $M_{nil}$ are obtained by Putcha, ~\cite[Theorem 3.1]{Putcha98}:
\begin{align}\label{A:interpretable}
M= \bigsqcup_{[\sigma]\in \mathcal{C}} X([\sigma]) \qquad \text{and}\qquad M_{nil}= \bigsqcup_{[\sigma]\in \mathcal{C}_{nil}} X([\sigma]),
\end{align}
where $X([\sigma])=\bigcup_{g\in G} g B \sigma B g^{-1}$.
The indexing set in the first decomposition, that is $\mathcal{C}$, is called the {\em Putcha poset of $M$}; 
it is defined as a certain subquotient of the Renner monoid of $M$. 
Here, the {\em Renner monoid} of $M$ is the finite inverse semigroup defined by $R:=\overline{N_G(T)}/T$, where $N_G(T)$ is the normalizer of a maximal torus $T$ that is contained in $B$, and the bar over it indicates the Zariski closure in $M$. 
Then $\mathcal{C}_{nil}$ is the set of nilpotent elements of $\mathcal{C}$.
The purpose of our article is to investigate these finite invariants for the ``dual canonical monoids.'' 
These monoids arise rather naturally as deformations of semisimple groups. 
Indeed, the {\em asymptotic semigroup of a semisimple group $G_0$}, denoted by $\textrm{As}(G_0)$,  
is the algebraic semigroup whose coordinate ring is the associated graded ring $\textrm{gr}\ \C[G_0]$,
where $\C[G_0]$ is the coordinate ring of $G_0$. 
The grading on $\C[G_0]$ is the one that comes from a well-known decomposition of $\C[G_0]$ as a $G_0\times G_0$-module. 
More precisely, we have $\C[G_0] = \bigoplus_{\chi \in \OO^+} V(\chi)\otimes V^*(\chi)$,
where $\OO^+$ is the semigroup of dominant weights, and $V(\chi)$ is the finite dimensional
irreducible representation of $G_0$ corresponding to the highest weight $\chi \in \OO^+$, and $V^*(\chi)$ is its dual. This remarkable algebraic semigroup is introduced by Vinberg in~\cite{Vinberg1,Vinberg2},
and studied by Rittatore~\cite{Rittatore:Thesis, Rittatore2001} from a viewpoint of spherical varieties.
By~\cite[Theorem 2]{Vinberg2}, we know that the union $M:=\textrm{As}(G_0) \sqcup G$,
where $G\cong \C^* \cdot G_0$, has the structure of a normal irreducible algebraic semigroup. 
An alternative construction of $M$ via one-parameter monoids is outlined in~\cite[Section 6.2]{Renner}. 
We note that since $G$ is present in it as a unit group, $M$ is in fact a semisimple monoid. 
For this reason, sometimes we refer to $M$ as the {\em asymptotic monoid of $G_0$}.
As we alluded before, $M$ is a dual canonical monoid. 
We will properly introduce the dual canonical monoids in the preliminaries section. 
\medskip

Let $W$ denote the Weyl group, $N_G(T)/T$. 
The {\em cross-section lattice of $(M,T)$} is the finite lattice of idempotents from $\overline{T}$, denoted $\Lambda$, 
such that $M= \bigsqcup_{e\in \Lambda} GeG$. 
For an idempotent $e\in \Lambda$, 
let $\mathcal{C}(e)$ denote the subposet of $\mathcal{C}$ defined by $\mathcal{C}(e):= \mathcal{C}\cap WeW$.
We are now ready to state our first main result.

\begin{Theorem}\label{T:firstmain-moregeneral}
Let $M$ be a dual canonical monoid.
If $e_\emptyset$ is the unique maximal element of $\Lambda \setminus \{1\}$, 
then $\mathcal{C}(e_\emptyset)$ is isomorphic to the opposite of the Bruhat-Chevalley order on $W$. 
Furthermore, under this isomorphism, the dimension of the corresponding subvariety, that is $\dim \overline{X([e_\emptyset w])}$,  
is given by $\dim G_0 - \ell(w)$, where $\ell(w)$ is the length of $w$ as an element of $W$. 
\end{Theorem}

We want to mention the fact that the first part of our Theorem~\ref{T:firstmain-moregeneral} follows easily from some general results that are proved by Putcha and Therkelsen. 
The real thrust of our result is its second assertion. 
A theorem of Putcha~\cite[Theorem 3.2 (ii)]{Putcha98} shows that the distinct irreducible components of the nilpotent variety of the dual canonical monoid are in one-to-one correspondence with the Coxeter elements of $W$. 
As a corollary of our theorem we obtain the following statement which strengthens Putcha's theorem. 

\begin{Corollary}\label{T:firstmain}
Let $M$ be a dual canonical monoid with unit group $G$.
Let $(W,S)$ denote the Coxeter system for the Weyl group of $G$. 
If $G_0$ denotes the semisimple part of $G$, then $M_{nil}$ is an equidimensional variety of dimension $\dim G_0 - |S|$. 
\end{Corollary}

In our next result, we will focus on the dual canonical monoid  of type $\textrm{A}_{n+1}$.
Then the unit group of $M$ is given by the (complex) general linear group $\mathbf{GL}_n$. 
In this case, the relevant combinatorics becomes especially concrete.
The Weyl group of $\mathbf{GL}_n$ is the symmetric group $S_n$. 
Let $R_n$ denote the {\em rook monoid} of $n\times n$ 0/1 matrices with at most one 1 in each row and each column. 
The rook monoid is the Renner monoid of the reductive monoid of $n\times n$ matrices, see~\cite{Renner86}. 
The Bruhat-Chevalley-Renner order, denoted $\leq$, is defined by the inclusion relations among the Zariski closures of 
$B\times B$-orbitsin $G$. 
We establish a connection between the Putcha poset $\mathcal{C}(e)$ of $M$  
and the Bruhat-Chevalley-Renner order on $R_n$.

\begin{Theorem}\label{T:main2}
Let $\mathcal{C}$ denote the Putcha poset of the dual canonical monoid whose unit group is $\mathbf{GL}_n$.
Let $k$ be a number such that $\lfloor n/2 \rfloor \leq k \leq n-1$. 
Let $(W,S)$ denote the Coxeter system of $\mathbf{GL}_n$, where $W$ is the symmetric group $S_n$, 
and $S:=\{s_1,\dots, s_{n-1}\}$ is the set of simple reflections that generate $W$. 
Let $e_I$ (resp. $W_I$) denote the idempotent determined by the set $I:=\{s_1,\dots, s_{k} \}$ (resp. the parabolic subgroup generated by $I$).
Then the following posets are isomorphic:
\begin{enumerate}
\item the opposite of the poset $W_I \backslash W / W_I$;
\item the Putcha subposet $\mathcal{C}(e_I)$;
\item $R_{\lfloor n/2 \rfloor -k}$.
\end{enumerate}
Here, $W_I \backslash W / W_I$ is the set of all two-sided cosets of $W_I$ in $W$. 
It is equipped with the order that is induced from the Bruhat-Chevalley order. 
\end{Theorem}
Various conjugacy actions on Renner monoids are investigated by Li, Li, and Cao in~\cite{LiLiCao2013}. 
In Section 4.1 of this reference, the authors show how to embed a Renner monoid into a rook monoid. 
It would be very interesting to use their result to extend our Theorem~\ref{T:main2} to the dual canonical monoids of other types. 
\medskip

Another goal of our paper is to initiate the study of the two-sided weak order, denoted by $\leq_{LR}$, 
on reductive monoids. 
We define our order by using the double Richardson-Springer monoid action on the Renner monoid $R$. 
This action respects the decomposition $R= \bigsqcup_{e\in \Lambda} WeW$.
When it is restricted to $W$, our two-sided weak order agrees with the ordinary two-sided weak order on $W$ viewed as a Coxeter group. 
It is easy to see from a simple example 
that the two-sided weak order on a Coxeter group is not a lattice. 
None the less, we show that for the dual canonical monoids, if $e$ is from $\Lambda \setminus \{1\}$, then 
$(WeW,\leq_{LR})$ is a lattice. 
Furthermore, we show that it is a distribute lattice if and only if $(WeW,\leq_{LR})\cong (WeW,\leq)$. 
\medskip

An important notion that is closely related to the geometry of weak order is the ``degree'' of a covering relation. 
Roughly speaking, it measures the generic degree of a morphism that is canonically attached to a covering relation in the weak order. 
This number (the degree) can be 0,1, or a positive power of 2. 
In this particle, we prove the following relevant theorem. 

\begin{Theorem}
Let $M$ be a dual canonical monoid, and let $W$ and $\Lambda$ denote, as before, 
the Weyl group and the cross-section lattice of $M$, respectively. If $e$ is an idempotent from 
$\Lambda \setminus \{1\}$, then all covering relations in $(WeW,\leq_{LR})$ have degree 1. 
\end{Theorem}

The two-sided weak order on the symmetric group $S_{n+1}$ is interesting by itself.
It turns out that there are many degree 2 covering relations in this case.

\begin{Theorem}
Let $W$ denote the symmetric group $S_{n+1}$. 
Then we have 
\begin{enumerate}
\item[(1)] the total number of covering relations in $(W,\leq_{LR})$ is $n^2 n!$;
\item[(2)] the number of covering relations of degree 2 in $(W,\leq_{LR})$ is $n n!$.
\end{enumerate}
\end{Theorem}

We are now ready to describe the individual sections of our paper. 
In the next preliminaries section we collect some well-known facts about the 
reductive monoids, Bruhat-Chevalley-Renner order, Putcha posets, and about nilpotent varieties.
The purpose of Section~\ref{S:Type} is to streamline some important structural results regarding 
the type map and the $G\times G$-orbits for a dual canonical monoid.
In Section~\ref{S:Rook}, we prove that the rook monoid appears as 
an interval in the Putcha poset of the dual canonical monoid with unit group $\mathbf{GL}_n$. 
In Section~\ref{S:Nilpotent}, we show that the nilpotent variety of the dual canonical monoid 
is equidimensional. In particular, we find the precise descriptions of certain intervals of $\mathcal{C}_{nil}$. 
In addition, we present a practical method (Theorem~\ref{T:practical}) 
for comparing two elements from different subposets $\mathcal{C}(e)$ and $\mathcal{C}(f)$.
The purpose of Section~\ref{S:RichardonSpringer} is to introduce the two-sided weak order on $WeW$. 
Also in this section, for dual canonical monoids, we present formulae for the cardinalities of the Renner monoid 
and of its set of idempotents (Corollary~\ref{C:formulae}).

\section{Preliminaries}\label{S:Preliminaries}

Let $G$ be a connected reductive group, let $T$ be a maximal torus, and let $B$ be a Borel subgroup of $G$ 
such that $T\subset B$. As before, let $W$ denote the Weyl group $N_G(T)/T$. 
The Bruhat-Chevalley order on $W$ is defined by $v\leq w \iff B\dot{v} B\subseteq \overline{B\dot{w} B}$,
where $\dot{v}$ and $\dot{w}$, respectively, are two elements from $N_G(T)$ representing the cosets $v$ and $w$. 
The bar on $B\dot{w}B$ indicates the Zariski closure in $G$. 
In the sequel, if a confusion is unlikely, then we will omit writing the dots on the representatives of the cosets.

For the poset $(W,\leq)$, the data of $(G,B,T)$ determines a Coxeter generating system $S$ 
and a length function $\ell : W \to \Z$, where, for $w\in W$, $\ell(w)$ is equal to the minimal number of 
simple reflections $s_{i_1},\dots, s_{i_r}$ from $S$ with $w=s_{i_1}\cdots s_{i_r}$.
A subgroup that is generated by a subset $I\subset S$ is denoted by $W_I$; 
it is called a {\it parabolic subgroup} of $W$. 
For $I\subseteq S$, we will denote by $D_I$ the following set:
\begin{align}
D_I:= \{ x\in W:\ \ell(xw ) = \ell(x) + \ell(w) \text{ for all } w\in W_I \}.
\end{align}

Let $M$ be a reductive algebraic group. This means that the unit group of $M$,
denoted by $G$, is a connected reductive algebraic group. 
Let $T$ be a maximal torus in $G$, and let $B$ be a Borel subgroup such that $T\subset B$. 
The following decompositions are well-known: 
\begin{enumerate}
\item $M = \bigsqcup_{r\in R} B\dot{r} B$ \qquad (the {\it Renner decomposition of $M$});
\item $M = \bigsqcup_{e\in \Lambda} G e G$ \qquad (the {\it Putcha decomposition of $M$}).
\end{enumerate}
In the first item, the parametrizing object $R$ is called the {\em Renner monoid of $M$}.
It is defined as the quotient, $R:= \overline{N_G(T)}/T$, where the bar on $N_G(T)$ indicates the Zariski closure in $M$. 
The dot on an element $r\in R$ indicates that we are taking its coset representative from $\overline{N_G(T)}$. 
The Renner monoid is a finite inverse semigroup with unit group $W$. 
A useful survey of Renner monoids can be found in~\cite{LiLiCao2014}.

The parametrizing object of the Putcha decomposition, that is $\Lambda$, is called the {\em cross-section lattice}
(or, the {\em Putcha lattice}) of $M$.
If $M$ has a zero, then $\Lambda$ can be defined as 
\begin{align*}
\Lambda := \{ e\in E(\overline{T}) :\ Be = eB e \},
\end{align*}
where $E(\overline{T})$ denotes the semigroup of idempotents of $\overline{T}$.
In fact, $\Lambda$ and $B$ determine each other, see~\cite[Theorem 9.10]{Putcha}.
This means also that the cross section lattice determines (and determined by) 
the set of Coxeter generators for $W$.

The set that is described in the next lemma is first used by Renner in~\cite{Renner86}, 
where, among other things, the Gauss-Jordan elimination method is generalized to arbitrary reductive monoids.
\begin{Lemma}
If $GJ=GJ(R,B)$ denotes the set $GJ:= \{ x \in R:\ Bx \subseteq x B \}$, then 
$GJ$ is a submonoid of $R$.
\end{Lemma}
\begin{proof}
Clearly, the neutral element of $R$ is contained in $GJ$.
If $x$ and $y$ are two elements from $GJ$,
then $B xy \subseteq xBy$ and $xBy \subseteq x y B$. It follows that $xy \in GJ$. 
\end{proof}

We will call $GJ$ the {\it Gauss-Jordan monoid} of $M$. 
Strictly speaking, $GJ$ is determined by $(T,B,M)$. 
Note that the unit group $W$ acts on $R$ by left multiplication, and $W\times W$ acts on $R$ by $(a,b)\cdot x = axb^{-1}$,
where $a,b\in W$ and $x\in R$. Then the $W$-orbits (resp. the $W\times W$-orbits) are parametrized by $GJ$ (resp. 
by $\Lambda$). Indeed, it is easy to see from~\cite[Proposition 8.9]{Renner} that 
\begin{align*}
| Wx \cap GJ | =1 \ \text{ for every $x\in R$.}
\end{align*}

The cross section lattice $\Lambda$ has a natural semigroup theoretic partial order:
\begin{align}\label{A:orderoncsl}
e \leq f \iff e = fe =ef\ \text{ for } e,f\in \Lambda.
\end{align}
We note that the order (\ref{A:orderoncsl}) agrees with the {\em Bruhat-Chevalley-Renner order on $R$},
which is defined by
\begin{align}\label{A:BCR}
x \leq y \iff Bx B \subseteq \overline{ByB}\ \text{ for } x,y\in R.
\end{align}

For $e\in \Lambda$, we have the following subgroups of $W$: 
\begin{enumerate}
\item $W(e):= \{ a\in W:\ ae=ea \}$,
\item $W^*(e):= \cap_{f\geq e} W(f)$,
\item $W_*(e):= \cap_{f\leq e} W(f) = \{ a\in W:\ ae= ea = e\}$.
\end{enumerate}
Then we know from~\cite[Chapter 10]{Putcha} that 
$W(e),W^*(e)$, and $W_*(e)$ are parabolic subgroups of $W$.
We know also that $W(e) \cong W^*(e) \times W_*(e)$. 
If $W(e)$ and $W_*(e)$ are parabolic subgroups of the form $W(e) = W_I$ and $W_*(e) = W_K$ for some subsets $I,K\subseteq S$, 
then we will use the following notation:
\begin{align*}
D(e) := D_I\qquad \text{and}\qquad D_*(e):= D_K.
\end{align*}

Let $B(S)$ denote the Boolean lattice of all subsets of $S$. 
The {\em type map} of $\Lambda$ is an order preserving map $\lambda : \Lambda \to B(S)$;  
it plays the role of a Coxeter-Dynkin diagram for $M$.
It is defined as follows. Let $e\in \Lambda$. Then $\lambda (e) :=\{s \in S: e s = s e \}$. 
Associated with $\lambda(e)$ are the following sets:
\begin{align*}
\lambda_*(e) := \cap_{f\leq e} \lambda (f)\qquad \text{and} \qquad \lambda^*(e):= \cap_{f\geq e} \lambda(f).
\end{align*}
Then we have 
\begin{align*}
W(e)= W_{\lambda(e)},\qquad W_*(e) = W_{\lambda_*(e)},\qquad W^*(e) = W_{\lambda^*(e)}.
\end{align*}

\textbf{Theorem/Definition (Pennell-Putcha-Renner):}
For each $x\in WeW$ there exist elements $a\in D_*(e), b\in D(e)$, which are uniquely determined by $x$, such that 
\begin{align}\label{A:std}
x= a e b^{-1}.
\end{align}
The decomposition of $x$ in (\ref{A:std}) is called the {\em standard form of $x$}.
Let $e, f$ be two elements from $\Lambda$. It is proven in~\cite{PPR97} that 
if $x= a eb^{-1}$ and $y= cf d^{-1}$ are two elements in standard form in $R$, then 
\begin{align}\label{A:BCR}
x \leq y \iff e\leq f,\ a \leq cw,\ w^{-1} d^{-1} \leq b^{-1} \qquad \text{for some $w\in W(f)W(e)$.}
\end{align}

Let $D(e)^{-1}$ denote the set $\{ b^{-1} :\ b\in D(e)\}$.
In this notation, the Gauss-Jordan monoid of $R$ has the following decomposition: 
\begin{align}\label{A:decomposition of GJ}
GJ = \bigsqcup_{e\in \Lambda} e D(e)^{-1}.
\end{align}
For $e,f\in \Lambda$, let $x$ (resp. $y$) be an element from $D(e)^{-1}$ (resp. from $D(f)^{-1}$). 
Then (\ref{A:BCR}) translates to the following statement:
\begin{align}\label{A:BCRonGJ}
ex \leq f y \iff y \leq w x \qquad \text{for some $w\in W(e)$}.
\end{align}

Another useful method for studying Bruhat-Chevalley-Renner order is introduced by Putcha in~\cite{Putcha04}.
Let $e$ and $f$ be two elements from $\Lambda$ such that $e\leq f$. 
Then Putcha defines the associated ``upward projection map'' $p_{e,f} : WeW\to WfW$.
He shows in that article that 
\begin{align*}
\sigma \leq \sigma ' \iff p_{e,f}(\sigma) \leq \sigma'\qquad \text{for $\sigma\in WeW$ and $\sigma' \in WfW$.}
\end{align*}
In the sequel, we will use the adaptation of these maps to the Putcha posets of dual canonical monoids.
This adaptation is already used by Therkelsen in~\cite{Therkelsen_Fields}.
In~\cite{Putcha04}, Putcha proved the following properties of the upward projection maps:
\medskip
Let $e,f\in \Lambda$ be such that $e\leq f$. Then
\begin{enumerate}
\item $p_{e,f} : WeW \to WfW$ is order preserving and $\sigma \leq p_{e,f}(\sigma)$ for all $\sigma \in WeW$.
\item If $\sigma \in WeW$, $\theta \in WfW$, then $\sigma \leq \theta \iff p_{e,f}(\sigma) \leq \theta$.
\item If $h\in \Lambda$ with $e\leq h \leq f$, then $p_{e,f}= p_{h,f} \circ p_{e,h}$.
\item $p_{e,f}$ is onto if and only if $\lambda_*(e)\subseteq \lambda_*(f)$.
\item $p_{e,f}$ is 1-1 if and only if $\lambda(f) \subseteq \lambda(e)$.
\end{enumerate}

\subsection{The conjugacy decomposition.}

The results that we mention in this subsection are obtained by Putcha in 
a series of papers,~\cite{Putcha87,Putcha98,Putcha05,Putcha08}.  

We maintain our notation from the previous subsection but let $M$ denote a reductive monoid with zero. 
It is easy to check that the relation $\sim$ defined by 
\begin{align}\label{A:conjugacy order}
ey \sim e'y' \iff wey w^{-1} = e'y' \ \text{ for some } w \in W
\end{align}
is an equivalence relation on $GJ$. 
Note that, if $ey\sim e'y'$, then we have $e=e'$.

\begin{Definition}
The set of equivalence classes of $\sim$ together with the order in (\ref{A:conjugacy order}) is called the {\em Putcha poset of $M$}, denoted by $\mathcal{C}$. 
For $e\in \Lambda$, we denote by $\mathcal{C}(e)$ the subposet $\mathcal{C}(e):=\{ [e v] :\ [ev] \in \mathcal{C} \}$.
We denote by $\mathcal{C}_{nil}$ the subposet consisting of nilpotent elements, 
\begin{align*}
\mathcal{C}_{nil}:= \{ [ey]\in \mathcal{C}:\ (ey)^k = 0 \text{ for some } k \in \N\}.
\end{align*}
Let us denote by $\mathcal{C}_{nil}(e)$ the subposet $\mathcal{C}_{nil}\cap \mathcal{C}(e)$.
\end{Definition}

The {\em conjugacy decomposition of $M$} is given by
\begin{align*}
M=\bigsqcup_{[ey] \in \mathcal{C}} X([ey]),
\qquad \text{where}\qquad X([ey]):= \bigcup_{g\in G} g Bey B g^{-1}.
\end{align*}
Note that we can define $X([\cdot ])$ not just for the elements of $GJ$ but for every element $\sigma \in R$ by the same definition, 
$X([\sigma]) := \bigcup_{g\in G} g B\sigma B g^{-1}$. 
Furthermore, it is easy to see that 
\begin{align*}
\tau \leq \sigma\implies X([\tau]) \subseteq \overline{X([\sigma])} \qquad \text{for every $\tau,\sigma \in R$}.
\end{align*}
Following Putcha we now define a partial order on $\mathcal{C}$:
\begin{align}\label{A:Putchaorder1}
[ey] \leq [e'y'] \iff X( [ey]) \subseteq \overline{X([e'y'])} \qquad \text{ for $[ey],[e'y']\in \mathcal{C}$}.
\end{align}
Then we have  
\begin{align*}
\overline{X([e'y'])}= \bigcup_{[ey]\leq [e'y']}  X([ey]).
\end{align*}
It turns out that the order (\ref{A:Putchaorder1}) is equivalent to the following partial order:
\begin{align}\label{A:Putchaorder}
[ey] \leq [e'y'] \iff wey w^{-1} \leq e'y'\qquad \text{for some $w\in W$.}
\end{align}
This is proved by Putcha in~\cite[Theorem 2.8]{Putcha05}.

\medskip

Recall that the {\em nilpotent variety} of $M$, denoted $M_{nil}$, is defined by 
\begin{align*}
M_{nil}:= \{ a\in M:\ a^k = 0 \text{ for some } k \in \N\}.
\end{align*}
The conjugacy decomposition of $M$ yields the following conjugacy decomposition of $M_{nil}$: 
\begin{align*}
M_{nil} = \bigsqcup_{[ev]\in \mathcal{C}_{nil}} X([ev]).
\end{align*}

A reductive monoid $M$ is called {\em $J$-coirreducible} if $\Lambda \setminus \{1\}$ has 
a unique maximal element, denoted by $e_{max}$. In this case, the {\em type of $M$} is defined 
as the subset $I:=\lambda(e_{max})$ in $S$.
A reductive monoid $M$ with a zero is called  {\em $J$-irreducible} if $\Lambda \setminus \{0\}$ 
has a unique minimal element, denoted by $e_{min}$. In this case, the {\em type of $M$} is defined 
as the subset $I:=\lambda(e_{min})$ in $S$.
The following theorem of Putcha will be useful in the sequel.

\begin{Theorem}\cite[Theorem 6.1]{Putcha08}\label{T:6.1}
Let $M$ be a $J$-coirreducible monoid of type $I$. Then
\begin{enumerate}
\item $M$ is semisimple, that is, the center of $G$ is one-dimensional; 
\item $e,e'\in \Lambda \setminus \{1\}$, then $e\leq e'$ if and only if $\lambda_*(e')\subseteq \lambda_*(e)$;
\item $e'\in \Lambda \setminus \{1\}$, then $\lambda^*(e) = \{ s \in I:\ s s' = s' s \text{ for every }s' \in \lambda_*(e) \}$;
\item If $K \subseteq S$, then $K=\lambda_*(e)$ for some $e\in \Lambda \setminus \{1\}$ 
if and only if no connected component of $K$ is contained in $I$; 
\item If $e\in \Lambda \setminus \{1\}$, then $|\lambda_*(e)| = crk (e) -1 = |S|- rk(e)$.
In particular, if $e\in \Lambda_{min}$, then $\lambda(e) =\lambda_*(e) = S\setminus \{s\}$ for some $s \in S$.
\end{enumerate}
\end{Theorem}

\begin{Example}
Let us denote by $\textbf{M}_n$ the monoid of $n\times n$ matrices. 
The unit group of $\textbf{M}_n$ is equal to $\mathbf{GL}_n$. 
It is well-known from linear algebra that the $\mathbf{GL}_n\times \mathbf{GL}_n$-orbits 
in $\textbf{M}_n$ are parametrized by the ranks of the matrices that are contained in the orbits.
It is also well-known that the Zariski closure of the orbit of a rank $k$ matrix in $\textbf{M}_n$ 
contains all other $n\times n$ matrices of lower ranks. Therefore, the cross-section
lattice of $\textbf{M}_n$ forms a chain of length $n+1$. In particular, we see that $\textbf{M}_n$ 
is a $J$-irreducible as well as a $J$-coirreducible monoid. 
To describe the corresponding types, let $S$ denote the set of simple transpositions of the set $\{1,\dots, n\}$. 
In other words, $S=\{s_1,\dots, s_{n-1}\}$, where $s_i$ ($i\in \{1,\dots, n-1\}$) is the permutation 
that interchanges $i$ and $i+1$, and $s_i(j)=j$ for $j\in \{1,\dots, n\}\setminus \{i,i+1\}$.  
Then $S$ is a Coxeter generating set for the Weyl group $W$ of $\mathbf{GL}_n$, which is isomorphic to the symmetric 
group of $\{1,\dots, n\}$. 
In this notation, the type of $\textbf{M}_n$, as a $J$-irreducible monoid, is given by 
$S\setminus \{s_1\}$. If we view $\textbf{M}_n$ as a $J$-coirreducible monoid, then its type is given by 
the set $S\setminus \{s_{n-1}\}$. 
\end{Example}

\begin{Definition}
Let $(W,S)$ be a Coxeter system. Let $s_1,\dots, s_n$ denote the elements of $S$. 
An element $v\in W$ is called {\em linear} if it is of the form $v:=s_{i_1}\cdots s_{i_p}$,
where $s_{i_1},\dots, s_{i_p}$ are all different from each other. 
A linear element is called a {\em Coxeter element} if $p= |S|$. 
\end{Definition}

In~\cite[Theorem 6.2]{Putcha08} Putcha shows that,
if $M$ is a $J$-coirreducible monoid of type $I$, then the distinct irreducible components of $M_{nil}$ are given by 
the Zariski closures $\overline{X([e_{max} x])}$, where $x$ is a Coxeter element of $W$ in $D_I^{-1}$.

\begin{Definition}
Let $M$ be a $J$-coirreducible monoid of type $I$.
If $I=\emptyset$, then $M$ is called a {\em dual canonical monoid}..
This means that $\lambda(e_{max}) = \emptyset$. 
In this case, we will denote $e_{max}$ by $e_\emptyset$. 
A canonical monoid is defined similarly; let $M$ be a $J$-irreducible monoid of type $I$. 
If $I=\emptyset$, then $M$ is called a {\em canonical monoid}.
\end{Definition}


\subsection{Double cosets.}\label{SS:double cosets}

Let $(W,S)$ be a Coxeter system, let $I$ and $J$ be two subsets from $S$. 
For $w\in W$, we denote by $[w]$ the double coset $W_I w W_J$. 
Let $\pi : W \rightarrow W_I\backslash W / W_J$
denote the canonical projection onto the set of 
$(W_I,W_J)$-double cosets. 
It turns out that the preimage in $W$ of every double coset 
in $W_I\backslash W / W_J$ is an interval with 
respect to Bruhat-Chevalley order, 
hence it has a unique maximal and a 
unique minimal element, see~\cite{Curtis85}.
Moreover, if $[w], [w']\in W_I\backslash W / W_J$ 
are two double cosets, $w_1$ and $w_2$ are the maximal 
elements of $[w]$ and $[w']$, respectively, 
then 
$w \leq w'$ if and only if $w_1 \leq w_2$,
see~\cite{HohlwegSkandera}.
Therefore, $W_I\backslash W / W_J$ has 
a natural combinatorial partial ordering 
defined by 
$
[w] \leq [w'] \iff w \leq w' \iff w_1 \leq w_2 
$
where $[w],[w'] \in W_I\backslash W / W_J$ and $w_1$ and 
$w_2$ are the maximal 
elements, $w_1 \in [w]$ and $w_2 \in [w']$.

Now let $[w]$ be a double coset from $W_I\backslash W / W_J$
represented by an element $w\in W$ such that 
$\ell(w) \leq \ell(v)$ for every $v\in [w]$. 
It turns out that the set of all such minimal length double coset representatives
is given by $D_I^{-1} \cap D_J$, the intersection of the set of minimal length left coset representatives of $W_I$ in $W$ and the set of minimal length right coset representatives of $W_J$ in $W$.
We will denote this intersection by $X_{I,J}^-$. 
Set $H:= I \cap w J w^{-1}$. 
Then $ uw \in D_J$ for $u\in W_I$ 
if and only if $u$ is a minimal length coset 
representative for $W_I/W_H$. 
In particular, every element of $W_I w W_J$ 
has a unique expression of the form 
$uwv$ with $u\in W_I$ is a minimal length 
coset representative of $W_I/W_H$, $v\in W_J$ and 
$\ell(uwv) = \ell(u)+\ell(w)+\ell(v)$. 

Another characterization of the sets $X_{I,J}^-$ is as follows. 
For $w\in W$, the {\em right ascent set} is defined as 
$
\mt{Asc}_R(w) = \{ s\in S :\ \ell(ws) > \ell(w) \}.
$
The {\em right descent set}, $\mt{Des}_R(w)$ is the complement $S \setminus \mt{Asc}_R(w)$. 
Similarly, the {\em left ascent set} of $w$ is 
\hbox{$\mt{Asc}_L(w) = \{ s\in S :\ \ell(sw) > \ell(w) \}$},
which is equal to $\mt{Asc}_R(w^{-1})$. 
Then we have 
\begin{align*}
X_{I,J}^- &= \{ w\in W:\ I\subseteq \mt{Asc}_L(w)\ \text{ and } J\subseteq \mt{Asc}_R(w) \}\\
&=  \{ w\in W:\ I^c\supseteq \mt{Des}_R(w^{-1})\ \text{ and } J^c\supseteq \mt{Des}_R(w) \}
\end{align*}
Let us point out that, in general, the Bruhat-Chevalley order on $X_{I,J}^-$ is a nongraded poset. 
For some special choices of $I$ and $J$, in type A, we determined the corresponding posets explicitly,
see~\cite{Can18,CanTien}.

\section{The Type Map of a Dual Canonical Monoid}\label{S:Type}

Most of the results in this section are well-known to the experts.
In fact, as observed by Therkelsen in~\cite{Therkelsen_Thesis}, 
the proofs of many of these results follow by duality from the corresponding facts 
that hold true in the canonical monoid case.
However, since they are important for our purposes, we provide direct proofs for completeness.

The Boolean lattice $B_n$ is the poset of all subsets of an $n$-element set which is ordered with respect to
the inclusions of subsets. 
The {\em opposite-Boolean lattice} is the opposite of the poset $(B_n,\subseteq)$. 
We will denote it by $B_n^{op}$. 
For $A,B\in B_n^{op}$, we have $A \leq B \iff A\supseteq B$.
To ease our notation, we denote the set $\{1,\dots, n\}$ by $[n]$. 

\begin{Lemma}\label{L:anti Boolean}
Let $P$ be a graded sublattice of $B_n^{op}$ with $\emptyset \in P$ and 
$[n]\in P$.
If for every element $I$ in $B_n^{op}$ there is a collection of elements $A_1,\dots, A_r$ in $P$ such that $\cap_{i=1}^r A_i = I$,
then $P=B_n^{op}$.
\end{Lemma}

\begin{proof}
Clearly our claim is true for $n=1$ as well as for $n=2$. We will prove the general case by induction,
so we assume that our lemma holds true for the opposite-Boolean poset $B_{n-1}^{op}$.

Now, let $P$ be a graded sublattice of $B_n^{op}$ which satisfies the hypothesis of our lemma.
Clearly, for every $i\in [n]$, the set $A_i:=[n] \setminus \{i\}$ is an element of $P$.
These are precisely the atoms in $P$. Note that if $K$ is a subset in $[n]$,
then $K= \cap_{i\in K} A_i$. 

Let $B(i)$ denote the opposite-Boolean sublattice in $B_n^{op}$ which consists of all 
subsets containing the element $i$. Then $A_1,\dots,A_{i-1},A_{i+1},\dots,A_n$ are 
elements of $B(i)$, and furthermore, any other element in $B(i)$ can be written as their intersections. 
Therefore, by our induction hypothesis the sublattice generated by 
$A_1,\dots,A_{i-1},A_{i+1},\dots,A_n$ is equal to $B(i)$. 
This arguments is true for all $i\in [n]$. 
Finally, we note that $\{\emptyset \} \cup \bigcup_{i\in [n]} B(i) = B_n^{op}$. This finishes the proof. 
\end{proof}

The opposite-Boolean lattice of subsets of $S$ will be denoted by $B^{op}(S)$.
Let $\Lambda$ be the cross-section lattice of a dual canonical monoid $M$.
When we want to be very precise, for $I \in B^{op}(S)$ such that $\lambda(e) = I$, we will write $e_I$ to specify $e$.

\begin{Proposition}\label{P:above}
Let $M$ be a dual canonical monoid. Then $\Lambda \setminus \{1\}$ is isomorphic to the 
opposite-Boolean lattice, $B^{op}(S)$.
\end{Proposition}

\begin{proof}
The cross section lattice of $M$ contains 0 as an element.  
It corresponds to $e_S$. Indeed, 
by part 4 of Theorem~\ref{T:6.1}, for $f\in \Lambda_{min}$, we have $\lambda(f) = S \setminus \{s \}$ for some $s \in S$.
This implies that $\lambda( 0 ) = S$.

Since $M$ is of type $\emptyset$, by part 3 of Theorem~\ref{T:6.1}, for any $K\subseteq S$ we have an idempotent 
$e\in \Lambda \setminus \{1\}$ such that $\lambda_*(e) = K$. 
We know that the type map $\lambda: \Lambda \to  B^{op}(S)$ is 1-1 in our case, therefore, 
$\Lambda \setminus \{1\}$ isomorphic to its image under $\lambda$.  
Since for every $e\in \Lambda \setminus \{1\}$, we have $\lambda_*(e) = \cap_{f\leq e} \lambda(f)$, 
we see that $\Lambda\setminus \{1\}$ satisfies the hypothesis of Lemma~\ref{L:anti Boolean}.
This finishes the proof. 
\end{proof}

\begin{Corollary}\label{C:*}
Let $M$ be a dual canonical monoid. 
Then $\lambda_*(e) = \lambda(e)$ for all $e\in \Lambda\setminus \{1\}$. 
Consequently, we have $W(e) = W_*(e)$ for all $e\in \Lambda\setminus \{1\}$. 
\end{Corollary}
\begin{proof}
Let $e$ be an idempotent in $\Lambda \setminus \{1\}$.
It follows from Proposition~\ref{P:above} that
if $f\in \Lambda \setminus \{1\}$ is such that $f\leq e$, then 
$\lambda(f) \supseteq \lambda(e)$. Therefore, $\lambda_*(e) = \cap_{f\leq e} \lambda(f) = \lambda(e)$. 
Our second assertion follows from the definitions of $W(e)$ and $W_*(e)$.
\end{proof}

For an element $e\in \Lambda$, let us denote by $P(e)$ and $P(e)^-$ the subgroups
\begin{align*}
P(e) = \{ g\in G:\ ge = ege \} \ \text{   and } \ P(e)^- = \{g\in G:\ eg = ege \}.
\end{align*}
Then $P(e)$ and $P(e)^-$ are opposite parabolic subgroups in $G$. 
The centralizer of $e$ in $G$ will be denoted by $C_G(e)$. In other words, 
we have $C_G(e):= \{ g\in G:\ ge=eg\}= P(e) \cap P(e)^-$.

\begin{Theorem}
Let $M$ be a dual canonical monoid, and let $e$ be an idempotent from 
the cross section lattice $\Lambda$ of $M$. 
If $B$ denotes the Borel subgroup that determines $\Lambda$, then the $G\times G$-orbit $GeG$ is a fiber bundle over 
$G/P(e)\times G/P(e)^-$ with fiber $eBe$ at the identity double coset $idP(e)\times idP(e)^-$.
\end{Theorem}

\begin{proof}
The following fibre bundle structure on $GeG$ is well-known:
\begin{align}\label{A:fibre}
eC_G(e) \to GeG \to G/P(e) \times G/P(e)^-.
\end{align}
A proof of it can be found in~\cite[Lemma 3.5 and 3.6]{CanRenner}.
Note that the second map in (\ref{A:fibre}) is given by $geh \mapsto (gP(e), h^{-1}P(e)^-)$ for $geh \in GeG$. 
By Corollary~\ref{C:*}, we know that $W(e) = W_*(e)=\{ w\in W:\ we=ew=e\}$. 
We know from~\cite[Proposition 10.9 (i)]{Putcha} that the Weyl group of $C_G(e)$ is given by $W(e)$.
Let $B_1$ denote the Borel subgroup of $C_G(e)$ such that $C_G(e) = B_1 W(e) B_1$
(the Bruhat-Chevalley decomposition for $C_G(e)$). Then we see that 
\begin{align*}
eC_G(e) = B_1 e W(e) B_1= B_1e W_*(e)B_1 =B_1eB_1=eB_1.
\end{align*}
But $eB_1 = e C_B(e)=eBe$ by~\cite[Corollary 7.2]{Putcha}. This finishes the proof.
\end{proof}

\begin{Corollary}\label{C:fiberT_0}
If $e$ is the idempotent $e= e_\emptyset$ in $\Lambda$, then 
$GeG$ is a torus fiber bundle over $G/B\times G/B^-$. 
More precisely, we have 
\begin{align*}
T_0 \to Ge_\emptyset G \to G/B\times G/B^-,
\end{align*}
where $T_0$ is the maximal torus of the derived subgroup of 
the unit group $G$. 
\end{Corollary}

\begin{proof}
This follows from the fact that if $e= e_\emptyset$, then $P(e)= B$, $P(e)^-=B^-$,
and \hbox{$C_G(e) = T$}. Finally, we note that $e_\emptyset T \cong T_0$ since 
$e_\emptyset$ is the maximal element of $\Lambda \setminus \{1\}$, and 
the height of $\Lambda \setminus \{1\}$ is equal to $\dim T_0$.
\end{proof}

\section{The Rook Monoid As an Interval}\label{S:Rook}

The following useful combinatorial result is first recorded by Therkelsen in his PhD thesis
\cite[Theorem 5.2.2]{Therkelsen_Thesis}.

\begin{Lemma}\label{L:Ryan}
Let $M$ be a dual canonical monoid with cross-section lattice $\Lambda$.
If $e$ is an element from $\Lambda\setminus \{1\}$, then 
$\mathcal{C}(e)$ is isomorphic to the dual of $W(e) \backslash W / W(e)$.
In other words, we have 
\begin{align*}
[ey] \leq [ex] \iff W(e) x W(e) \leq W(e) y W(e) \iff x \leq y,
\end{align*}
for $x,y\in D_*(e) = D(e)\cap D(e)^{-1}$.
\end{Lemma}
Here, it is a natural question to ask for which idempotents $e\in \Lambda \setminus \{1\}$
the double coset $W(e) \backslash W / W(e)$ is graded. 
For $e= e_\emptyset$ this is the case. 
We will reprove this result in the proof of our Theorem~\ref{T:firstmain-moregeneral}.
In type A, our results in~\cite{Can18} shows that if $e=e_{S\setminus \{s\}}$, then 
$W(e)\backslash W / W(e)$ is a graded lattice. We anticipate this result will hold true in other 
types as well.

The {\em rook monoid} on the set $\{1,\dots, n\}$, denoted by $R_n$, is the full inverse semigroup 
of injective partial transformations $\{1,\dots, n\} \to \{1,\dots, n\}$. It is the Renner monoid 
of the reductive monoid of $n\times n$ matrices. The unit group of $R_n$ is the symmetric 
group $S_n$. Let $w$ be a permutation from $S_{n}$. The {\em one-line notation} for $w$
is a string of numbers $w_1\dots w_{n}$, where $w_i = w(i)$ for $i\in \{1,\dots, n\}$. 
In a similar manner, the {\em one-line notation} for $\sigma \in R_n$ is a string 
of numbers $\sigma_1\dots \sigma_n$, where, for $i\in \{1,\dots, n\}$, $\sigma_i = \sigma(i)$ 
if $\sigma(i)$ is defined; otherwise $\sigma_i= 0$. 
For example, $\sigma = 02501$ is the injective partial transformation
$\sigma : \{2,3,5\}\to \{1,2,3,4,5\}$ with $\sigma(2) = 2$, $\sigma(3)=5$, and $\sigma(5)=1$.

Let $\sigma = \sigma_1\dots \sigma_n$ and $\tau= \tau_1\dots \tau_n$ be two elements from 
$R_n$. 
We will write $\widetilde{\sigma_i}$ for the non-increasing rearrangement of the string
$\sigma_1\sigma_2\dots \sigma_i$. For example, if $\sigma = 02501$, then 
$\widetilde{\sigma_4}= 5200$. 
If $a:= a_1\dots a_m$ and $b:=b_1\dots b_m$ are two strings of integers of the same length, 
then we will write $a \leq_c b$ if $a_i \leq b_i$ for all $i\in \{1,\dots, m\}$. 
The following characterization of the Bruhat-Chevalley-Renner order is proven 
in~\cite{CanRenner12}:
\begin{align}\label{A:Deodhar}
\tau \leq \sigma \iff \widetilde{\tau_i} \leq_c \widetilde{\sigma_i} \ \text{ for all }i\in \{1,\dots, n\}.
\end{align}
Our next result describes a surprising connection between $R_n$ and the Putcha monoid 
of the dual canonical monoid with unit group $\mathbf{GL}_n$.

\begin{Theorem}\label{T:rook}
Let $W$ denote the symmetric group $W= S_{2m}$. 
If $I$ denotes the subset $\{s_1,\dots, s_{m-1} \}$ in $S=\{s_1,\dots, s_{2m-1}\}$,
then the opposite of the poset $W_I \backslash W / W_I$, or equivalently,
the Putcha subposet $\mathcal{C}(e_I)$ is isomorphic to 
the poset $(R_m,\leq)$.
\end{Theorem}

\begin{proof}
First, we will determine the elements of $D_*(e_I)$. 
Let $w=w_1\dots w_{2m}$ be an element from $D_*(e_I)$.
Notice that the set $I$ indicates the positions of the descents in $w$;
if $s_i\in I$, then $w_i > w_{i+1}$. 
Since $w^{-1}$ is also in $D(e_I)$, we see that
if $w_{i_1}=2m, w_{i_2}=2m-1,\dots, w_{i_m}=m+1$, then $i_1 < \cdots < i_m$. 
At the same time, $w$ is of minimal possible length. These requirements
imply that the intersection 
$
\{1,\dots, m\} \cap \{i_1,\dots, i_m\} = \{i_1,\dots, i_k\}
$ 
uniquely determines $w$; we place $2m,\dots, 2m-k+1$ at the positions $i_1,\dots, i_k$,
and we place $2m-k+2,\dots, m+1$ at the positions $m+1,m+2,\dots, 2m-k$.
The numbers $i_1,\dots, i_k$ are placed, in a decreasing order, at the positions
$2m-k+1,\dots, 2m$. The remaining entries are filled in the increasing order with what remains 
of $1,2,3,\dots,2m$.
But now such a permutation, $w\in S_{2m}$ defines a unique partial permutation 
with its first $m$ entries; we define
$\sigma=\sigma(w)$ by $\sigma_i:=w_i -i$ for $i\in \{1,\dots, m\}$. 
It is not difficult to show conversely that any $\sigma \in R_m$ gives 
a permutation $w=w(\sigma) \in D_*(e_I)\subset S_{2m}$.
Furthermore, it is now clear from (\ref{A:Deodhar}) that,  
for two elements $\tau$ and $\sigma$ from $R_m$, 
$\tau \leq \sigma$ if and only if $w(\tau) \leq w(\sigma)$. 
This finishes the proof.
\end{proof}

The proofs of the next two corollaries follow from the proof of Theorem~\ref{T:rook}.

\begin{Corollary}\label{C:rook}
Let $W$ denote the symmetric group $W= S_{2m+1}$. 
If $I$ denotes the subset $\{s_1,\dots, s_{m} \}$ in $S=\{s_1,\dots, s_{2m}\}$,
then the opposite of the poset $W_I \backslash W / W_I$, or equivalently,
the Putcha subposet $\mathcal{C}(e_I)$ is isomorphic to 
the poset $(R_m,\leq)$.
\end{Corollary}

\begin{Corollary}
Let $W$ denote the symmetric group $W= S_n$. 
Let $k$ be a number such that $\lfloor n/2 \rfloor \leq k \leq n-1$.
If $I$ is the subset $\{s_1,\dots, s_{k} \}$ in $S=\{s_1,\dots, s_{n-1}\}$, 
then the opposite of the poset $W_I \backslash W / W_I$, or equivalently,
the Putcha subposet $\mathcal{C}(e_I)$ is isomorphic to $(R_{\lfloor n/2 \rfloor -k},\leq)$.
\end{Corollary}

\begin{figure}[htp]
\begin{center}
\scalebox{.65}{
\begin{tikzpicture}[scale=.34]

\node at (0,-5) (a0) {$14523$};
\node at (0,0) (a) {$12543$};
\node at (-5,5) (b1) {$12453$};
\node at (5,5) (b2) {$12534$};
\node at (-5,10) (c1) {$12354$};
\node at (5,10) (c2) {$12435$};
\node at (0,15) (d) {$12345$};

\draw[-, thick] (a0) to (a);
\draw[-, thick] (a) to (b1);
\draw[-, thick] (a) to (b2);
\draw[-, thick] (b1) to (c1);
\draw[-, thick] (b1) to (c2);
\draw[-, thick] (b2) to (c2);
\draw[-, thick] (b2) to (c1);
\draw[-, thick] (c1) to (d);
\draw[-, thick] (c2) to (d);

\end{tikzpicture}
}
\caption{The Putcha subposet $\mathcal{C}(e_{\{s_1,s_2\}})$ is isomorphic to 
the rook monoid $(R_2,\leq)$.}\label{F:1}
\end{center}
\end{figure}

\section{The Nilpotent Variety of a Dual Canonical Monoid}\label{S:Nilpotent}

Let $M$ be a dual canonical monoid, and let $\mathcal{C}$ denote the corresponding Putcha monoid. 
Let $[ev]$ ($v\in D(e)^{-1}$) be an element from $\mathcal{C}$. 
Putcha proved in~\cite[Theorem 4.2]{Putcha08} that $[ev]\in \mathcal{C}_{nil}$ if and only if $\textrm{supp}(v) \nsubseteq \lambda(f)$ for all 
 $f\in \Lambda_{min}$ with $f\leq e$.  
 Also, we know from the previous section that for such $f$, $\lambda(f) = S \setminus \{s\}$ for some 
$s \in S$, and $f \leq e$ if and only if $\lambda(e) \subseteq \lambda(f)$. 
Therefore, $\textrm{supp}(v)$ contains every $s$ that lies in the complement of the set $\lambda(e)$.
In other words, we have
\begin{align}\label{A:supp v}
\textrm{supp}(v) \supseteq S\setminus \lambda(e).
\end{align}
As a consequence of this observation, 
we identify the maximal elements of the subposet $\mathcal{C}_{nil}(e) \subseteq \mathcal{C}(e)$ 
for $e\in \Lambda \setminus \{1\}$.

\begin{Proposition}\label{P:commutes}
Let $K$ be a subset of $S$. Then the set of maximal elements of the poset $\mathcal{C}_{nil}(e_K)$
consists of  elements of the form $[e_K s_{i_1}s_{i_2}\cdots s_{i_k}]$, where $\{s_{i_1},\dots, s_{i_k} \} = S\setminus K$.
In particular, $\mathcal{C}_{nil}(e_K)$ has a unique maximal element if and only if 
$s_i s_j = s_j s_i$ for all $s_i,s_j$ in $S\setminus K$.
\end{Proposition}

\begin{proof}
Let $[ex]$ and $[ey]$ be two elements from $\mathcal{C}_{nil}(e_K)$.
By Lemma~\ref{L:Ryan}, $[ex]\leq [ey]$ if and only if $y\leq x$. 
Therefore, by (\ref{A:supp v}), the maximal elements of $\mathcal{C}_{nil}(e_K)$ are of the form 
$[ey]=[e (s_{i_1}s_{i_2}\cdots s_{i_k})]$, where $\{i_1,\dots, i_k\} = S\setminus \lambda (e_K)$.
The second claim is obvious.
\end{proof}

\begin{Corollary}
Let $e_K$ be a minimal nonzero idempotent from $\Lambda \setminus \{1\}$. 
Then $\mathcal{C}_{nil}(e_K)$ has a unique maximal and a unique minimal element.
\end{Corollary}

\begin{proof}
If $e_K$ is a minimal nonzero element in $\Lambda\setminus \{1\}$, then by Proposition~\ref{P:commutes}
we know that $K=S \setminus \{ s \}$ for some $s \in S$. Therefore, $S\setminus K  = \{s\}$.
In other words, $\mathcal{C}_{nil}(e_K)$ has a unique maximal and a unique minimal element.
\end{proof}

\begin{Remark}
In type $A$, for $K=S\setminus \{s\}$, the poset $\mathcal{C}(e_K)$, hence $\mathcal{C}_{nil}(e_K)$, is a chain.
This statement holds true in some other types as well, see~\cite[Proposition 3.2]{Littelmann} and~\cite[Theorem 2.3]{Stembridge04}.
\end{Remark}

Let $e_I$ and $e_J$ be two different elements from $\Lambda\setminus \{1\}$. 
Comparisons between the elements belonging to $\mathcal{C}(e_I)$ and $\mathcal{C}(e_J)$ are described by another result of Therkelsen. By using Therkelsen's and Putcha's results, we observe that the lower interval $[e_S, e_\emptyset w_0]$ is a Boolean lattice. 

\begin{Proposition}
The interval between $[e_\emptyset w_0]$ and $[e_S]$ ($e_S=0$) in
$\mathcal{C}$, hence in $\mathcal{C}_{nil}$, is isomorphic to $B^{op}(S)$.  
\end{Proposition}

\begin{proof}
Let $I$ be a subset of $S$, and let $[e_I y]$ be the minimal element of interval $\mathcal{C}(e_I)$. 
Then $[e_I y]\in \mathcal{C}_{nil}$. Let $J$ be another subset of $S$. 
If $[e_J z]$ is the minimal element of $\mathcal{C}(e_J)$, then we will prove that 
$
J\subseteq I \iff [e_I y] \leq [e_J z].
$
Clearly, $(\Leftarrow)$ is true. 
To prove the other direction, we will prove the stronger statement that $e_I y \leq e_J z$ in the Bruhat-Chevalley-Renner
order. By~\cite[Lemma 2.1 (i)]{Putcha05} this will show that $ [e_I y] \leq [e_J z]$ in $\mathcal{C}$. 
To prove the latter statement, first, we will show that 
\begin{align}\label{A:showthis}
p_{e_I, e_\emptyset} (e_I y) = e_\emptyset w_0.
\end{align} 
By the last part of Theorem~\ref{T:6.1} and Corollary~\ref{C:*}, we know that the upward projection maps are one-to-one. 
Thus, we conclude that $ [e_I y] \leq [e_J z]$ in $\mathcal{C}$. 
Now (\ref{A:showthis}) can be seen directly from the description of the Bruhat-Chevalley-Renner order (\ref{A:BCR}) as follows. 
We write $w_0$ in the form $w^{-1}y^{-1}=w_0$ for some $w^{-1} \in W(e_I)$. 
Then (\ref{A:BCR}) shows that $[e_I y ] \leq [e_\emptyset w_0]$, hence, it shows that (\ref{A:showthis}).
This finishes the proof.
\end{proof}

We now proceed to prove our Theorem~\ref{T:firstmain-moregeneral}. 
Let us recall its statement for convenience: 

Let $M$ be a dual canonical monoid.
If $e_\emptyset$ is the unique maximal element of $\Lambda \setminus \{1\}$, 
then $\mathcal{C}(e_\emptyset)$ is isomorphic to the opposite of the Bruhat-Chevalley order on $W$. 
Furthermore, for every $w\in W$, the dimension of the corresponding Putcha sheet, that is, $\dim X([e_\emptyset w])$, 
is given by $\dim G_0 - \ell(w)$, where $\ell(w)$ is the length of $w$ as an element of $W$.

\begin{proof}[Proof of Theorem~\ref{T:firstmain-moregeneral}]
For the idempotent $e_\emptyset$, we have the following identifications: 
\begin{align*}
W(e_\emptyset) = W^*(e_\emptyset) = W_*(e_\emptyset) = \{id\}  
\end{align*}
and 
\begin{align*}
D(e_\emptyset)^{-1} = D_*(e_\emptyset) = W.
\end{align*}
Since $W(e_\emptyset)=\{id\}$, by~\cite[Theorem 2.2]{Putcha05}, for $y,z\in D(e_\emptyset^{-1})$, the following conditions are equivalent: 
\begin{enumerate}
\item[(i)] $[e_\emptyset y ] \leq [e_\emptyset z ]$ 
\item[(ii)] $z \leq y$.
\end{enumerate} 
In particular, since $\leq$ is a partial order, the equality $[e_\emptyset y ] = [e_\emptyset z ]$ holds if and only if the equality $y=z$ holds. 
It follows that $\mathcal{C}(e_\emptyset)$ is isomorphic to the opposite of the Bruhat-Chevalley order on $W$, denoted $W^{op}$. 
This finishes the proof of our first assertion.

The maximal element of $\mathcal{C}(e_\emptyset)$ is then $[e_\emptyset]$, which corresponds to the open stratum 
$X([e_\emptyset]) = \bigcup_{g\in G} g B e_\emptyset B g^{-1}$.
Therefore, the dimension of $X([e_\emptyset])$ is given by 
\begin{align*}
\dim X([e_\emptyset]) = \dim G e_\emptyset G = \dim M - 1 = \dim G_0.
\end{align*}

We assume that the length of $y$ as an element of $W$ is $r:= \ell(y)$. 
Let $s_{i_1} s_{i_2} \cdots s_{i_r}$ be a reduced expression for an element $y\in W$. 
Then we have an associated increasing saturated chain in $W^{op}$, 
\begin{align*}
y=s_{i_1}s_{i_2}\cdots s_{i_{r}} \lessdot s_{i_1}s_{i_2}\cdots s_{i_{r-1}} \lessdot \cdots \lessdot s_{i_1} s_{i_2} \lessdot s_{i_1} \lessdot id.
\end{align*}
Corresponding to this chain we have a chain of varieties, 
\begin{align*}
\overline{ X([e_\emptyset y]) } \subsetneq \overline{ X([e_\emptyset s_{i_1}s_{i_2}\cdots s_{i_{r-1}]}) } \subsetneq \cdots \subsetneq 
\overline{ X([e_\emptyset s_{i_1}s_{i_2}]) } \subsetneq \overline{ X([e_\emptyset s_{i_1}]) } \subsetneq \overline{ X([e_\emptyset])} = \overline{G e_\emptyset G}.
\end{align*}
It follows that 
\begin{align}\label{A:upperbound}
\dim \overline{ X([e_\emptyset y])} \leq \dim G_0 - r.
\end{align}
Thus, it remains to show that $\dim \overline{ X([e_\emptyset y])}\geq \dim G_0 - r$.
We will consider the orbit $G\cdot ey = \{ ge\dot{y} g^{-1} :\ g\in G\}$ in $X([e_\emptyset y])$

Let $\varphi : GeG\to G/B\times G/B^-$ denote the fibration that is in Corollary~\ref{C:fiberT_0}; it is given by 
\begin{align*}
\varphi (aeb) = (aB, bB^-) \qquad (a,b\in G).
\end{align*}
Thus, for $g\in G$, we have $\varphi (g\cdot e\dot{y}) = (gB, \dot{y}g^{-1}B^-)$.
By using Chevalley's big cell theorem~\cite[Proposition 28.5]{Humphreys}, if we write $g$ in the form $g= ub$, where $u\in U^-$ and $b\in B$,
then we see that 
\begin{align*}
\varphi (g e\dot{y} g^{-1}) = (gB, g\dot{y}^{-1}B^-) = (uB,ub \dot{y}^{-1} B^-).
\end{align*} 
Notice that if we fix $u$ in $U^-$ and vary $b\in B$, then we obtain the subset 
$\{uB\}\times u B\dot{y}^{-1} B^-$ in $G/B\times G/B^-$. 
Clearly, this subset is isomorphic via the second projection $G/B\times G/B^-\to G/B^-$ to the $u$-translate of the opposite Schubert cell $B\dot{y}^{-1} B^-$ in $G/B^-$. 
Indeed, the following identifications in $G/B^-$ are easily checked:
\begin{align}\label{A:easilychecked}
B\dot{y}^{-1} B^- = w_0 w_0 B w_0 w_0 \dot{y}^{-1} B^- = w_0 B^- w_0\dot{y}^{-1} B^-
\end{align}
Therefore, the dimension of the opposite Schubert cell $B\dot{y}^{-1} B^-$ is given by 
\begin{align*}
\dim B^- w_0\dot{y}^{-1} B^- &= \ell(w_0) - \ell(y^{-1})  \\
&= \ell(w_0) - \ell(y)\\
&= \ell(w_0) - r.
\end{align*}
Next, we let $u$ vary in $U^-$ to get the subset $U^-B\times U^-B\dot{y}^{-1} B^-$ in $G/B\times G/B^-$. 
Let $p_1 : G/B\times G/B^- \to G/B$ be the first projection. 
Clearly, for every $u\in U^-$, the dimension of the intersection $p_1^{-1}(uB) \cap (U^-B\times U^-B\dot{y}^{-1} B^-)$
is $\ell(w_0) - \ell(y)$. In other words, we have a constant fiber dimension. 
Finally, since we have $p_1(U^-B\times U^-B\dot{y}^{-1} B^-) = U^-B$ is open in $G/B$, we see that the dimension of the image of $p_1$ is $\ell(w_0)$. 
It follows that the dimension of $\varphi( X([e_\emptyset y]))$ is given by $2\ell(w_0) - r$. 
At the same time, $\varphi$ has constant fiber dimension $\dim T_0$. 
Therefore, the dimension of $X([e_\emptyset y])$ is at least $2\ell(w_0) - r+ \dim T_0$. 
It follows that the dimension of the orbit $G\cdot e\dot{y}$ is $2\ell(w_0) -r + \dim T_0$.
Hence the dimension of $X([e_\emptyset y])$, is at least $2\ell(w_0) -r + \dim T_0$.
But it is easy to check that $\dim G_0 = 2\ell(w_0) + \dim T_0$. 
In other words, we have $\dim X([e_\emptyset y]) \geq \dim G_0 - r$. 
This, together with the upper bound in (\ref{A:upperbound}) finishes the proof of our second assertion, hence, the theorem is proved. 
\end{proof}

The proof of the following corollary follows easily from the fact that the length of a Coxeter element $y\in W$ is equal to 
number of simple generators of $W$. 
We omit its details. 

\begin{Corollary}
Let $M$ be a dual canonical monoid.
Let $S$ denote the set of Coxeter generators of $W$ relative to $B$. 
Then $M_{nil}$ is an equidimensional variety of dimension $\dim G_0 - |S|$. 
\end{Corollary}

We now present a useful criterion for comparing two elements from different subposets $\mathcal{C}(e)$ and $\mathcal{C}(f)$.
A similar result is obtained by Therkelsen in his thesis,~\cite[Theorem 6.2.2]{Therkelsen_Thesis}.

\begin{Theorem}\label{T:practical}
Let $[ex ]$ and $[fy]$ be two elements from $\mathcal{C}$,
where $e,f\in \Lambda\setminus \{1\}$, and $x \in D(e)^{-1}, y\in D(f)^{-1}$. 
Then 
$[ex ] \leq [f y]$ if and only if $e\leq f$ in $\Lambda$ and there exists 
$w\in W(e)$ such that 
if we write $x w$ in the form $v_1 v_2$ with $v_1 \in W(e)$ and $v_2\in D(e)^{-1}$, then $v_2 \geq y$.
In particular, for every $e,f\in \Lambda$ with $e\leq f$ and $v\in W(e)$, $x \in D(e)^{-1}$, we have 
$[ex]=[evx] \leq [fx]$. 
\end{Theorem}

\begin{proof}
Let $S$ denote as before the Coxeter generators of $W$. 
Let $I$ and $J$ denote the subsets from $S$ that correspond to $e$ and $f$, respectively. 
By (\ref{A:Putchaorder}), $[ex] \leq [f y]$ if and only if $w exw^{-1} \leq fy$ for some $w\in W$. 
Let us write $w$ in the form $w_1w_2$ with $w_2\in W(e)$ and $w_1\in D(e)$.
Since $W(e) = W_*(e)$ by Corollary~\ref{C:*}, we have $wex w^{-1} = w_1 w_2 ex w^{-1} = w_1 ex w^{-1}$. 
Likewise, we write $xw^{-1}$ in the form $v_1 v_2$, where $v_1\in W(e)$ and $v_2\in D(e)^{-1}$.
Then we have $w_1 ex w^{-1}=w_1 e v_2$. 
Recall that $D(e) = D_*(e)$. 
It follows from the Pennell-Putcha-Renner theorem (\ref{A:BCR}) that 
$w ex w^{-1} = w_1 e v_2 \leq fy$ if and only if there exists $u\in W(f)W(e)$ such that 
$e\leq f$, $w_1 \leq u$, $u^{-1} y \leq v_2$. The last inequality is equivalent to the inequality $y^{-1} u \leq v_2^{-1}$. 
Note that $e\leq f$ if and only if $I \supseteq J$ if and only if $W(e) \supseteq W(f)$ if and only if $D(e) \subseteq D(f)$.
Under these conditions, $y^{-1} u \leq v_2^{-1}$ holds only if $u=id=w_1$, and therefore, $w= w_2$. 
In other words, $[ex ] \leq [f y]$ if and only if there exists $w\in W(e)$ (hence $w^{-1}\in W(e)$) such that 
if we write $x w^{-1} = v_1 v_2$ with $v_1 \in W(e)$ $v_2\in D(e)^{-1}$, then $v_2 \geq y$.
This finishes the proof of our first assertion. 

For our second claim, we observe that, since $D(e) \subseteq D(f)$ implies that $x \in D(f)^{-1}$. 
Now by letting $w$ denote the identity element, and $v_2:=x$, we see that $x w^{-1} = v_2 = x$, 
hence, by the first part of the proposition, we obtain $[ex ] \leq [fx]$. 
Finally, since $W(e) = W_*(e)$, we have $ev = e$ for every $v\in W(e)$. This finishes the proof.  
\end{proof}

\begin{Example}
Let $G_0$ denote the the exceptional simple algebraic group of type $\textrm{G}_2$.
Then the Weyl group $W$ of $G_0$ is isomorphic to the dihedral group of order 12. 
To setup our notation, in Figure~\ref{F:Directed}, we depicted the directed Coxeter-Dynkin diagram of $W$. 
\begin{figure}[htp]
\begin{center}

\scalebox{.75}{
\begin{tikzpicture}
\pgfarrowsdeclare{arcs}{arcs}{...}
{
\pgfsetdash{}{0pt} 
\pgfsetroundjoin 
\pgfsetroundcap 
\pgfpathmoveto{\pgfpoint{-10pt}{10pt}}
\pgfpatharc{180}{270}{10pt}
\pgfpatharc{90}{180}{10pt}
\pgfusepathqstroke
}
\filldraw (-1,0) circle (.15cm);
\filldraw (1,0) circle (.15cm);
\node at (-1,0.35) (a') {$s$};
\node at (1,0.35) (b') {$t$};
\draw[-arcs,line width=1pt] (1,0) -- (-.25,0);
\draw[-, line width=1pt] (-.25,0) -- (-1,0);
\draw[-, line width=1pt] (1,-0.1) -- (-1,-0.1);
\draw[-, line width=1pt] (1,0.1) -- (-1,0.1);
\end{tikzpicture}
}
\caption{The directed Coxeter-Dynkin diagram of type $\textrm{G}_2$.}
\label{F:Directed}
\end{center}
\end{figure}
The simple reflection $s$ corresponds to the short simple root in $\textrm{G}_2$.
We have the relations $s^2=t^2=1$ and $(st)^6=1$. 
The longest element of $W$ is given by $w_0:= st st st = tststs$.
It is easy to read from this data that $G_0$ is 14 dimensional, 
and a maximal unipotent subgroup in $G_0$ is 6 dimensional.
In particular, the full flag variety of $G_0$ is 8 dimensional.  
Let $M$ denote the dual canonical monoid for $G_0$. 
Let $S$ denote the set $\{s,t\}$. The cross-section lattice of $M$ has four elements, 
which are parametrized by the subsets of $S$,
$\Lambda = \{ e_\emptyset, e_{\{ s \}}, e_{\{ t\}}, e_S \}$.
Since $W(e_\emptyset) = \{1\}$, $W(e_{\{s\}}) = W_{\{ s\}}$,$W(e_{\{t\}}) = W_{\{t\}}$,
and $W(e_S) = W$, we see that the corresponding minimal length left coset representatives are given by 
\begin{align*}
D(e_\emptyset) &= W \Rightarrow D(e_\emptyset)^{-1} = W,\\
D(e_{\{t\}}) &= \{ 1, s, ts, sts, tsts,ststs \} \Rightarrow D(e_{\{t\}})^{-1} = \{ 1, s, st, sts, stst,ststs \},\\
D(e_{\{s\}}) &= \{1,t, st, tst, stst, tstst \} \Rightarrow D(e_{\{s\}})^{-1} = \{1,t, ts, tst, tsts, tstst \},\\
D(e_S ) &= \{1\} \Rightarrow D(e_S)^{-1} = \{1\}.
\end{align*} 
By using this data and Proposition~\ref{T:practical}, 
we obtain the complete picture of the Putcha poset of the dual canonical monoid associated with 
the simple algebraic group $\textrm{G}_2$. 
In Figure~\ref{F:G2}, we depict $\mathcal{C}_{nil}$ and $\mathcal{C}$, 
where the former poset is shown in thicker fonts. 
Notice in this case that both of the posets $\mathcal{C}$ and $\mathcal{C}_{nil}$ are graded posets. 
Notice also that the interval $[ [e_\emptyset w_0], [e_\emptyset]]$ in $\mathcal{C}$ is isomorphic to the (opposite) Bruhat-Chevalley order on the Weyl group of $\textrm{G}_2$, and the interval $[ [e_S], [e_\emptyset w_0]]$ is the opposite of the Boolean lattice of $S$. 

\begin{figure}[htp]
\begin{center}

\scalebox{.76}{
\begin{tikzpicture}[scale=.25]
\node[blue] at (0,-15) (a0) {$[e_S]$};
\node[blue] at (-5,-10) (b0) {$[e_{\{s\}} tstst]$};
\node[blue] at (5,-10) (c0) {$[e_{\{t\}} ststs ]$};

\node[blue] at (0,-5) (a1) {$[e_{\emptyset}w_0]$};
\node[blue] at (-5,0) (a2) {$[e_{\emptyset}ststs]$};
\node[blue] at (5,0) (a3) {$[e_{\emptyset}tstst]$};
\node[blue] at (-5,5) (a5) {$[e_{\emptyset}tsts]$};
\node[blue] at (5,5) (a4) {$[e_{\emptyset}stst]$};
\node[blue] at (-5,10) (a6) {$[e_{\emptyset}sts]$};
\node[blue] at (5,10) (a7) {$[e_{\emptyset}tst]$};
\node[blue] at (-5,15) (a9) {$[e_{\emptyset}ts]$};
\node[blue] at (5,15) (a8) {$[e_{\emptyset}st]$};
\node at (-5,20) (a10) {$[e_{\emptyset}s]$};
\node at (5,20) (a11) {$[e_{\emptyset}t]$};
\node at (0,25) (a12) {$[e_{\emptyset}]$};

\node[blue] at (15,-5) (c1) {$[e_{\{t\}} stst]$};
\node[blue] at (15,0) (c2) {$[e_{\{t\}} sts]$};
\node[blue] at (15,5) (c3) {$[e_{\{t\}} st]$};
\node[blue] at (15,10) (c4) {$[e_{\{t\}} s]$};
\node at (15,15) (c5) {$[e_{\{t\}}]$};

\node[blue] at (-15,-5) (b1) {$[e_{\{s\}} tsts]$};
\node[blue] at (-15,0) (b2) {$[e_{\{s\}} tst]$};
\node[blue] at (-15,5) (b3) {$[e_{\{s\}} ts]$};
\node[blue] at (-15,10) (b4) {$[e_{\{s\}} t]$};
\node at (-15,15) (b5) {$[e_{\{s\}}]$};

\draw[-, ultra thick, blue] (b1) to (a2);
\draw[-, ultra thick, blue] (b2) to (a5);
\draw[-, ultra thick, blue] (b3) to (a6);
\draw[-, ultra thick, blue] (b4) to (a9);
\draw[-] (b5) to (a10);

\draw[-, ultra thick, blue] (c1) to (a3);
\draw[-, ultra thick, blue] (c2) to (a4);
\draw[-, ultra thick, blue] (c3) to (a7);
\draw[-, ultra thick, blue] (c4) to (a8);
\draw[-] (c5) to (a11);

\draw[-, ultra thick, blue] (b0) to (b1);
\draw[-, ultra thick, blue] (b1) to (b2);
\draw[-, ultra thick, blue] (b2) to (b3);
\draw[-, ultra thick, blue] (b3) to (b4);
\draw[-] (b4) to (b5);

\draw[-, ultra thick, blue] (c0) to (c1);
\draw[-, ultra thick, blue] (c1) to (c2);
\draw[-, ultra thick, blue] (c2) to (c3);
\draw[-, ultra thick, blue] (c3) to (c4);
\draw[-] (c4) to (c5);

\draw[-, ultra thick, blue] (a0) to (b0);
\draw[-, ultra thick, blue] (a0) to (c0);
\draw[-, ultra thick, blue] (b0) to (a1);
\draw[-, ultra thick, blue] (c0) to (a1);
\draw[-, ultra thick, blue] (a1) to (a2);
\draw[-, ultra thick, blue] (a1) to (a3);
\draw[-, ultra thick, blue] (a2) to (a4);
\draw[-, ultra thick, blue] (a2) to (a5);
\draw[-, ultra thick, blue] (a3) to (a4);
\draw[-, ultra thick, blue] (a3) to (a5);
\draw[-, ultra thick, blue] (a4) to (a6);
\draw[-, ultra thick, blue] (a4) to (a7);
\draw[-, ultra thick, blue] (a5) to (a6);
\draw[-, ultra thick, blue] (a5) to (a7);
\draw[-, ultra thick, blue] (a6) to (a8);
\draw[-, ultra thick, blue] (a6) to (a9);
\draw[-, ultra thick, blue] (a7) to (a8);
\draw[-, ultra thick, blue] (a7) to (a9);
\draw[-] (a8) to (a10);
\draw[-] (a8) to (a11);
\draw[-] (a9) to (a10);
\draw[-] (a9) to (a11);
\draw[-] (a10) to (a12);
\draw[-] (a11) to (a12);
\end{tikzpicture}
}
\caption{The Putcha poset of the dual canonical monoid for $\textrm{G}_2$.}
\label{F:G2}
\end{center}
\end{figure}

\end{Example}

\section{A Richardson-Springer Monoid Action}\label{S:RichardonSpringer}

Let $M$ be a dual canonical monoid, and let $M_{nil}$ denote its nilpotent variety. 
The irreducible components of $M_{nil}$ are indexed by the Coxeter elements of the 
Weyl group of the unit group of $M$. It is well-known that all Coxeter elements are conjugate to each other.
However, they (Coxeter elements) do not necessarily form a single conjugacy class in a Weyl group.
Therefore, the conjugation action of $W$ on the set of Coxeter elements does not give
an additional structure to study the geometries of $M$ and $M_{nil}$. 
For this reason, 
in this section, we look at the actions of a certain deformation of the group ring of a Weyl group
on the reductive monoids. 
The structure that we will use is given by a finite monoid that is canonically associated with 
$W$, which is first used by Richardson and Springer in~\cite{RS90} for studying the weak order 
on symmetric varieties. 

\begin{Definition}\label{D:RS monoid}
Let $(W,S)$ be a Coxeter group. 
The {\em Richardson-Springer monoid} $O(W)$ of $W$ is the quotient of the free monoid generated by $S$ modulo the relations $s^2 = s$ for $s \in S$ and
\begin{equation}\label{eqn:braid rels}
stst \cdots = tsts \cdots
\end{equation}
for $s, t \in S$, where both sides of (\ref{eqn:braid rels}) are the product of exactly order of $st$ many elements.
\end{Definition}

$O(W)$ is a finite monoid, and its elements are in canonical bijection with the elements of $W$.  
We write $m(w)$ for the element of $O(W)$ corresponding to $W$.
If $w = s_1 s_2 \cdots s_l$ is any reduced expression of $w \in W$, 
then $m(w) = m(s_1) m(s_2) \cdots m(s_l)$. Furthermore, 
for $s \in S$ and $w \in W$, we have
\begin{align}\label{A:action of O(W)}
m(s) m(w)=
\begin{cases} 
m(sw) & \text{ if $\ell(sw) > \ell(w)$;} \\
m(w) & \text{ if $\ell(sw) < \ell(w)$.}
\end{cases}
\end{align}
From now on, we write $w$ for $m(w)$ when discussing an element $w \in O(W)$.
There is a useful geometric interpretation of (\ref{A:action of O(W)}). 
Let $X$ be a $G$-variety, and let $B$ be a Borel subgroup in $G$.
The set of all nonempty, irreducible, $B$-stable subvarieties of $X$ will be denoted by 
$\mathcal{B}(B:X)$. 
For $w\in W$, let $X_w$ denote the Zariski closure of $BwB$ in $G$. 
Clearly, every closed irreducible $B\times B$-subvariety of $G$ is of this type. 
For $w,w'\in W$, we set $X_{w*w'} := X_w X_w'$.
It is not difficult to check that if $s\in S$, $w\in W$, then 
$
X_{s *w} = X_{m(s)m(w)},
$
and that $X{s*w} \neq X_w$ if and only if $\ell (sw) = \ell(w) +1$.

Next, we will introduce the Richardson-Springer monoid 
action on $\mathcal{B}(B:X)$. 
For $Y\in \mathcal{B}(B:X)$, 
we have a morphism defined by the action, $\pi: G\times Y \to X$ $(g,z)\mapsto gz$.
Let $w$ be an element from $O(W)$. The restriction of $\pi$ to $X_w \times Y$
is equivariant with respect to $B$-action that is given by $b\cdot (a,z) := (ab^{-1},bz)$
for $b\in B$ and $(a,z)\in X_w\times Y$. 
Passing to the quotient, we get a new morphism 
$
\pi_{Y,w} : X_w \times^B Y \longrightarrow \overline{X_w Y}. 
$
Following~\cite{Knop}, let us denote $\overline{X_w Y}$ by $w*Y$. 
Next definitions are due to Brion~\cite[Section 1]{Brion98_Infinity}.
Since $1\in X_w$, we always have $Y\subseteq w*Y$. Note that it may happen that 
$Y= w*Y$ although $w\neq 1$. 
Note also that since $X_w/B$ is a complete variety, $\pi_{Y,w}$ is a proper map, hence, it is surjective. 
If the morphism $\pi_{Y,w}$ is generically finite, then we will denote 
the degree of $\pi_{Y,w}$ by $\deg(Y,w)$; if it is not generically finite, then 
we set $\deg(Y,w) :=0$. 
Finally, we define the {\em $W$-set of $Y$}, denoted $W(Y)$, 
as the set of $w$ from $O(W)$ such that $\pi_{Y,w}$ is generically finite and 
$\overline{Bw Y}$ is $G$-invariant. 
The following facts are proven in~\cite[Lemma 1.1]{Brion98_Infinity} 
\begin{Lemma}
Let $Y$ be a variety from $\mathcal{B}(B:X)$.
\begin{enumerate}
\item For any $\tau,w\in W$ such that $\ell(w\tau )= \ell(w) + \ell(\tau)$, we have 
$
d(Y, \tau w) = d(Y, \tau ) d(\overline{BwY}, \tau ).  
$
\item For any $w\in W$ such that $\overline{BwY}$ contains only finitely many $B$-orbits 
the integer $d(Y,w)$ is either 0 or a power of 2. 
\item For any $w\in W$ such that $d(Y,w)\neq 0$, we have 
\begin{align*}
W( \overline{BwY})= \{ \tau \in W:\ \ell(\tau w ) = \ell(\tau ) + \ell(w) \text{ and } \tau w \in W(Y)\}.
\end{align*}
\item The set $W(Y)$ is nonempty.
\item Assume that $X=G/P$, where $P$ is a parabolic subgroup with $B\subset P$, 
and with a Levi subgroup $L$ such that $T\subset L$. 
If $Y=\overline{BwP}/P$ with $\tau$ is a minimal length coset representative for $W_L$ in $W$, 
then $W(Y) = \{ w_0 w_{0,L} w^{-1} \}$, where $w_{0,L}$ denotes the longest element of 
$W_L$. Moreover, we have $d(Y, w_0 w_{0,L} w^{-1}) = 1.$
\end{enumerate}
\end{Lemma}

\begin{Definition}\label{D:weak}
Let $Y_1$ and $Y_2$ be two elements from $\mathcal{B}(B:X)$. 
We will write 
\begin{align}\label{A:weak order}
\text{$Y_1 \leq Y_2$ if  $Y_2 = w*Y_1$ for some $w\in O(W)$.}
\end{align} 
From now on, we will refer to the partial order that is defined by 
the transitive closure of the relations in (\ref{A:weak order}) the {\em weak order on $X$}.
If $Y_2 = \overline{Bs Y_1}$ for some $s\in S$ and $Y_2 \neq Y_1$, 
then we will call the cardinality $|W(Y_2)|$, the {\em degree} of the covering relation $Y_1 < Y_2$. 
In this case, we will write $\deg(Y_1,Y_2)$ for $|W(Y_2)|$.
\end{Definition}

\begin{Example} 
Let $I$ be a subset of $S$, and let $P=B W_I B$ denote the corresponding 
parabolic subgroup in $G$. We set $X:= G/P$, and let $Y$ be a Schubert variety in $X$
such that $Y= \overline{B w P}/P$, where $w\in D_I$. 
For $s\in S$, either $\dim s*Y = \dim Y$ or $\dim s*Y = \dim Y +1$. 
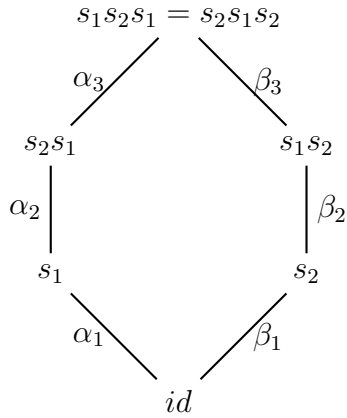
In the latter case, $\ell(sw) = \ell(w)+1$, and we get a covering relation for the left weak order on $D_I$.
In other words, the weak order on $X$ as defined in Definition~\ref{D:weak} agrees with 
the well-known left weak order on $D_I$.  
Furthermore, Brion's lemma shows that all covering relations in this case have degree 1. 
\end{Example}

Now we will apply this development in the setting of reductive monoids. 
By Bruhat-Chevalley-Renner order, we know that the set $\mathcal{B}(B\times B: M)$ 
is parametrized by the Renner monoid of $M$. 
Therefore, if we view $M$ as a $G\times G$-variety, 
then we have the ``doubled'' Richardson-Springer monoid action, 
\hbox{$*:O(W\times W) \times \mathcal{B}(B\times B: M) \to \mathcal{B}(B\times B:M)$,}
which is defined as follows: Let $s\in S$ and $\sigma \in R$. Then 
\begin{align}\label{A:action on R1}
(s,1)* \sigma =
\begin{cases} 
s \sigma & \text{ if $\ell(s\sigma ) > \ell(\sigma)$,} \\
\sigma & \text{ if $\ell(s\sigma ) \leq \ell(\sigma )$,}
\end{cases}
\end{align}
and 
\begin{align}\label{A:action on R2}
(1,s)* \sigma =
\begin{cases} 
\sigma s & \text{ if $\ell(\sigma s) > \ell(\sigma)$,} \\
\sigma & \text{ if $\ell(\sigma s ) \leq \ell(\sigma )$.}
\end{cases}
\end{align}
The operation in (\ref{A:action on R1}) corresponds to 
$Y \rightsquigarrow \overline{BsB Y}$, where $Y = \overline{B\sigma B}$,
and the operation in (\ref{A:action on R2}) corresponds 
to $Y \rightsquigarrow \overline{YBsB}$.
We will denote the weak order on $M$ by $(R,\leq_{LR})$.
This notation will be justified in the sequel.

Let $X$ be a $G$-variety, and let $Z$ be an element from $\mathcal{B}(B:X)$.
If $Z \subseteq Y$, where $Y$ is a $G$-orbit closure in $X$, then $w*Z \subseteq Y$
for all $w\in O(W)$. 
Consequently, we see that the weak order on $X$ is a disjoint union 
of various weak order posets, one for each $G$-orbit. 
It is easy to verify that 
$
(R,\leq_{LR}) = \bigsqcup_{e\in \Lambda} (WeW,\leq_{LR}).
$
Note that if $e$ is the neutral element of $G$, then we have 
$(WeW,\leq_{LR}) \cong (W,\leq_{LR})$. On the latter poset, the subscript $LR$
in the partial order stands for the two-sided weak order on the Coxeter group,
so, our choice of notation is consistent with the notation in the literature. 
We will denote the left (resp. right) weak order by $\leq_L$ (resp. by $\leq_R$).

Our next result shows the special nature of the Bruhat-Chevalley-Renner order on the dual canonical monoids. 

\begin{Proposition}\label{P:it is a lattice}
Let $\Lambda$ be a cross-section lattice of a reductive monoid, and 
let $e$ be an element from $\Lambda \setminus \{1\}$. 
If $\lambda^*(e) = \emptyset$, then we have
the following poset isomorphisms: 
\begin{enumerate}
\item[(1)] $(WeW,\leq) \cong (D(e),\leq)\times (D(e),\leq)^{op}$,  
\item[(2)] $(WeW,\leq_{LR}) \cong (D(e),\leq_L)\times (D(e),\leq_L)^{op}$.
\end{enumerate}
Furthermore, $(WeW,\leq_{LR})$ is a lattice.
\end{Proposition}
\begin{proof}
We start with the proof of (2). 
If $\lambda^*(e) = \emptyset$, then by using the standard forms of elements in $WeW$,
we see that $WeW = D(e) e D(e)^{-1}$. 
Let $\sigma = x e y$ and $\sigma ' = x' e y'$ be two elements from $D(e) e D(e)^{-1}$. 
Then $\sigma$ covers $\sigma'$ in $\leq_{LR}$ if and only if there exists $s\in S$
such that either $(s,1)*\sigma'= \sigma$, or $(1,s)* \sigma' = \sigma$. 
In the former case, $x$ covers $x'$ in $\leq_L$ and $y= y'$; 
in the latter case $y'$ covers $y$ in $\leq_R$,
hence $y'^{-1}$ covers $y^{-1}$ in $\leq_L$, and we have $x= x'$. 
This shows that the posets $(WeW,\leq_{LR})$ 
and $(D(e),\leq_L)\times (D(e),\leq_L)^{op}$ are canonically isomorphic. 
It is well known that the weak order on a quotient is a lattice. Since 
a product of two lattices is a lattice, the proof of (2) is finished. 

To prove (1), as before, let $\sigma = x e y$ and $\sigma ' = x' e y'$ be two elements from $D(e) e D(e)^{-1}$. 
By (\ref{A:BCR}) we know that $\sigma \leq \sigma'$ if and only if there exists $w\in W(e)$ such that 
$x \leq x' w$ and $w^{-1} y'^{-1} \leq y^{-1}$, or $y'w \leq y$. 
But since $x,x',y$, and $y'$ are from $D(e)$, and $w\in W(e)$, these inequalities simplify to give the 
inequalities $x \leq x'$ and $y'\leq y$. 
This finishes the proof. 
\end{proof}

As a consequence of Proposition~\ref{P:it is a lattice}, we see formulae for the cardinalities
of the Renner monoid of a dual canonical monoid and its idempotents. 
\begin{Corollary}\label{C:formulae}
Let $R$ denote the Renner monoid of a dual canonical monoid. 
Then the cardinality of $R$ is given by the formula,
$|R| = |W| + \sum_{I \subseteq S} \left( \frac{ |W| }{ |W_I |} \right)^2$,
where $W$ is the unit group of $R$, and $S$ is the set simple generators of $W$. 
The number of idempotents of $R$ is given by the summation $1 + \sum_{I \subseteq S} \left( \frac{ |W| }{ |W_I |} \right)$.
\end{Corollary}
\begin{proof}
The proof of our first claim follows from Proposition~\ref{P:it is a lattice} 
and the fact that $R = \bigsqcup_{e\in \Lambda} WeW$. 
For the second claim,we note that all idempotents of $R$ are of the form $w e w^{-1}$, where $w\in W$ and $e\in \Lambda$. 
Since $wew^{-1} = v e v^{-1}$, where $v$ is the minimal length coset representative for $wW_*(e)$,
we see that the number of conjugates of $e$ is equal to $|D(e)| = | W/ W_*(e)|$. 
The rest of the proof follows from this observation. 
\end{proof}

We note our formula in Corollary can be viewed as a special case of a theorem of Li, Li, and Cao~\cite[Theorem 2021]{LiLiCao2006}.

Let $W$ be an irreducible Coxeter group, and 
let $I$ be a subset of the set of simple roots $S$ for $W$. 
The set $D_I$ ($\cong W/W_I$) is said to be {\em minuscule} 
if the parabolic subgroup $W_I$ is the stabilizer of a ``minuscule'' weight.
Here, a weight $\nu$ is said to be {\em minuscule} if there is a representation of 
a semisimple linear algebraic group $G$ with Weyl group $W$ whose set of 
weights is the $W$-orbit of $\nu$.

The following result can be seen as an extension of~\cite[Theorem 7.1]{StembridgeFully}
into our setting.

\begin{Corollary}\label{C:distributive}
Let $e$ be an idempotent from a cross-section lattice of a reductive monoid $M$.
We assume that $e$ is not the neutral element. 
If $\lambda_*(e) \notin\{ \emptyset,S\}$ and $\lambda^*(e)= \emptyset$, then the following are equivalent.
\begin{enumerate}
\item $(WeW,\leq)$ is a lattice.
\item $(WeW,\leq)$ is a distributive lattice.
\item $(WeW,\leq_{LR})$ is a distributive lattice.
\item $(WeW,\leq_{LR}) = (WeW,\leq)$.
\item $D(e)$ is minuscule. 
\end{enumerate}
\end{Corollary}

\begin{proof}
Let $A$ and $B$ be two posets. The product poset $A\times B$ is a distributive lattice 
if and only if both of $A$ and $B$ are distributive lattices. 
Also, $A$ is a distributive lattice if and only if its opposite $A^{op}$ is a distributive lattice. 
Now, by Proposition~\ref{P:it is a lattice}, 
$(WeW,\leq_{LR})$ is always a lattice, and 
$(D(e),\leq)$ is a lattice if and only if $(WeW,\leq)$ is a lattice.
The rest of the proof follows from the proof of~\cite[Theorem 7.1]{StembridgeFully}.
\end{proof}

Next, we discuss the degrees of the covering relations for $\leq_{LR}$. 
Clearly, $(s,1)*1 = s = (1,s)*1$, therefore, the degree of the covering relation 
$1 < s$ in $(W,\leq_{LR})$ is always 2. 
\begin{Proposition}\label{P:degree 2}
Let $x,y$ be two elements from $W$. If $x$ is covered by $y$ in $(W,\leq_{LR})$,
then the degree of the covering relation is either 1 or 2. In the latter case,
there exist $s,s'\in S$ such that $y= (s,1)*x = (1,s')*x$. 
\end{Proposition}

\begin{proof}
Clearly, if $(s,1)*x = (s',1)*x = y$ for some $s,s'\in S$, then $s=s'$. 
Similarly, if $(1,s)*x = (1,s')*x = y$ for some $s,s'\in S$, then $s=s'$. 
Therefore, if the degree of $x<y$ is at least 2, then we can only have 
$(s,1)*x = (1,s')*x=y$ for some $s,s'\in S$. 
By the same argument, if they exist, then $s$ and $s'$ are unique. 
Therefore, the degree of a covering relation in $(W,\leq_{LR})$ is always $\leq 2$. 
\end{proof}

\begin{Theorem}\label{T:degree>1}
Let $\Lambda$ be a cross-section lattice of a reductive monoid,
and let $e$ be an element from $\Lambda \setminus \{1\}$. 
Then $\lambda^*(e) \neq \emptyset$ if and only if there is a covering relation 
$x<_{LR}y$ in $WeW$ such that $\deg(x,y) = 2$. 
\end{Theorem}
\begin{proof}
If $\lambda^*(e) \neq \emptyset$, then we know that $W^*(e) \neq \emptyset$,
hence, there is a simple reflection $s$ in $W^*(e)$ such that 
$es =se \neq e$. But this means that $\deg ( e,es) =2$.

Conversely, let $x$ be an element in $WeW$.
Let $a e b^{-1}$ be the standard form of $x$, 
where $a\in D_*(e)$ and $b\in D(e)$. 
By Proposition~\ref{P:degree 2}, if a covering relation $x<_{LR} y$ in $WeW$ has degree 2, then 
$(s,1)* x= (1,s')* x =y$ for some $s,s'\in S$. 
By the uniqueness of the standard form for the elements of $R$, 
the equality $sae b^{-1}  = a e b^{-1} s'$ implies that $s$ commutes 
with $a$ and $se=e$. Similarly, $s'$ commutes with $b^{-1}$ and 
$es'=e$. Since $R$ is a symmetric inverse semigroup, these equalities 
imply that $se=e=es$ and $es'=e=s'e$, hence $W^*(e)\neq \emptyset$. 
In other words, $\lambda^*(e)\neq \emptyset$. 
\end{proof}

To complete the hypothesis of Theorem~\ref{T:degree>1}, let $e$ denote the neural element, so, we have $WeW=W$.
It is easy to verify that in most Weyl groups there is a degree 2 covering relation.
For example, see Figure~\ref{F:S4}, where we depict $(S_4,\leq_{LR})$ together with all of its degree 2 covering relations. 
\begin{figure}[htp]
\begin{center}

\scalebox{.65}{
\begin{tikzpicture}[scale=.4]

\node at (0,0) (a) {$1234$};

\node at (-8,5) (b1) {$1243$};
\node at (0,5) (b2) {$1324$};
\node at (8,5) (b3) {$2134$};

\node at (-16,10) (c1) {$1342$};
\node at (-8,10) (c2) {$1423$};
\node at (0,10) (c3) {$2143$};
\node at (8,10) (c4) {$2314$};
\node at (16,10) (c5) {$3124$};

\node at (-20,15) (d1) {$1432$};
\node at (-12,15) (d2) {$2341$};
\node at (-4,15) (d3) {$2413$};
\node at (4,15) (d4) {$3142$};
\node at (12,15) (d5) {$3214$};
\node at (20,15) (d6) {$4123$};

\node at (-16,20) (e1) {$2431$};
\node at (-8,20) (e2) {$3241$};
\node at (0,20) (e3) {$3412$};
\node at (8,20) (e4) {$4132$};
\node at (16,20) (e5) {$4213$};

\node at (-8,25) (f1) {$3421$};
\node at (0,25) (f2) {$4231$};
\node at (8,25) (f3) {$4312$};

\node at (0,30) (g) {$4321$};

\draw[-, double, thick, blue] (a) to (b1);
\draw[-, double, thick, blue] (a) to (b2);
\draw[-, double, thick, blue] (a) to (b3);

\draw[-, thick] (b1) to (c1);
\draw[-, thick] (b1) to (c2);
\draw[-, double, thick, blue] (b1) to (c3);

\draw[-, thick] (b2) to (c1);
\draw[-, thick] (b2) to (c2);
\draw[-, thick] (b2) to (c4);
\draw[-, thick] (b2) to (c5);

\draw[-, double, thick, blue] (b3) to (c3);
\draw[-, thick] (b3) to (c4);
\draw[-, thick] (b3) to (c5);

\draw[-, double, thick, blue] (c1) to (d1);
\draw[-, thick] (c1) to (d2);
\draw[-, thick] (c1) to (d4);

\draw[-, double, thick, blue] (c2) to (d1);
\draw[-, thick] (c2) to (d3);
\draw[-, thick] (c2) to (d6);

\draw[-, thick] (c3) to (d3);
\draw[-, thick] (c3) to (d4);

\draw[-, thick] (c4) to (d2);
\draw[-, thick] (c4) to (d3);
\draw[-, double, thick, blue] (c4) to (d5);

\draw[-, thick] (c5) to (d4);
\draw[-, double, thick, blue] (c5) to (d5);
\draw[-, thick] (c5) to (d6);

\draw[-, thick] (d1) to (e1);
\draw[-, thick] (d1) to (e4);

\draw[-, double, thick, blue] (d2) to (e1);
\draw[-, double, thick, blue] (d2) to (e2);

\draw[-, thick] (d3) to (e1);
\draw[-, thick] (d3) to (e3);
\draw[-, thick] (d3) to (e5);

\draw[-, thick] (d4) to (e2);
\draw[-, thick] (d4) to (e3);
\draw[-, thick] (d4) to (e4);

\draw[-, thick] (d5) to (e2);
\draw[-, thick] (d5) to (e5);

\draw[-, double, thick, blue] (d6) to (e4);
\draw[-, double, thick, blue] (d6) to (e5);

\draw[-, thick] (e1) to (f1);
\draw[-, thick] (e1) to (f2);
\draw[-, thick] (e2) to (f1);
\draw[-, thick] (e2) to (f2);

\draw[-, double, thick, blue] (e3) to (f1);
\draw[-, double, thick, blue] (e3) to (f3);

\draw[-, thick] (e4) to (f2);
\draw[-, thick] (e4) to (f3);

\draw[-, thick] (e5) to (f2);
\draw[-, thick] (e5) to (f3);

\draw[-, double, thick, blue] (f1) to (g);
\draw[-, double, thick, blue] (f2) to (g);
\draw[-, double, thick, blue] (f3) to (g);

\node at (-4.5,3) {$(s_3,s_3)$};
\node at (0,2.5) {$(s_2,s_2)$};
\node at (4.5,3) {$(s_1,s_1)$};

\node at (-2.5,9.5) {$(s_1,s_1)$};
\node at (2.5,9.5) {$(s_3,s_3)$};

\node at (-18.5,12) {$(s_3,s_2)$}; 
\node at (-16,14) {$(s_3,s_3)$}; 

\node at (15,11.5) {$(s_2,s_1)$}; 
\node at (10,11.5) {$(s_3,s_2)$}; 

\node at (-15,18.5) {$(s_3,s_3)$}; 
\node at  (-9,18.5) {$(s_3,s_1)$}; 

\node at (16,16.5) {$(s_2,s_3)$}; 
\node at (18,18) {$(s_1,s_2)$}; 

\node at (-4,27.5) {$(s_3,s_1)$};
\node at (0,27.5) {$(s_2,s_2)$};
\node at (4,27.5) {$(s_1,s_3)$};

\node at (2.5,21.5) {$(s_3,s_1)$};
\node at (-2.5,21.5) {$(s_3,s_1)$};

\end{tikzpicture}
}
\caption{The two-sided weak order on $S_4$ and its double edges.}\label{F:S4}
\end{center}
\end{figure}

\begin{Corollary}\label{C:DCM is nice}
If $M$ is a dual canonical monoid and $e$ is an idempotent from 
$\Lambda \setminus \{1\}$, then all covering relations 
in $(WeW,\leq_{LR})= (D(e) e D(e)^{-1},\leq_{LR})$ are of degree 1. 
\end{Corollary}

\begin{proof}
This follows from Theorem~\ref{T:degree>1} and the fact that in a dual canonical
monoid we have $\lambda^*(e) = \emptyset$ for all $e\in \Lambda \setminus \{1\}$, see part 3 of Theorem~\ref{T:6.1}.
\end{proof}

In the rest of this section, we will consider the monoid $\textbf{M}_n$. 
Let $\textbf{B}_n$ denote the Borel subgroup consisting of upper triangular matrices in $\mathbf{GL}_n$. 
Then the corresponding cross-section lattice is given by $\Lambda := \{e_0=0, e_1,\dots, e_{n-1},e_n = 1\}$,
where $e_i$ ($i\in \{0,1,\dots,n\}$) is the diagonal matrix 
\begin{align*}
e_i:= \textrm{diag}(1,\dots, 1,0,\dots, 0)\ \text{ with $i$ 1's.}
\end{align*}
\begin{Proposition}
Let $e_i$ be a nonzero element from the cross-section lattice $\Lambda$ of $\textbf{M}_{n+1}$. 
Let $W$ denote $S_{n+1}$, the Weyl group of the unit group of $\textbf{M}_{n+1}$.
Then $i=1$ if and only if $\deg(x,y) =1$ for all covering 
relations $x<_{LR}y$ in $We_i W$. 
Furthermore, in this case, poset $(We_1 W , \leq_{LR})$ is isomorphic to 
$(W e_1 W, \leq)$.
\end{Proposition}

\begin{proof}
For the monoid $\textbf{M}_{n+1}$, 
it is easy to check that $\lambda^*(e_i) \neq \emptyset$ if and only if \hbox{$i\in \{2,\dots, n+1\}$.}
It is also easy to check that $\lambda_*(e_1) = \{s_2,\dots, s_n\}$. 
Therefore, our first claim follows from Theorem~\ref{T:degree>1},
and our second claim follows from Corollary~\ref{C:distributive}. 
\end{proof}

Let $x=x_1\dots x_{n+1}$ be a permutation in one-line notation.
A {\em right ascent} in $x$ is 
a string of two consecutive integers $\alpha:= i\, i+1$ such that $x_{i+1} > x_i$. 
A {\em small (right) ascent} in $x$ is 
a string of two consecutive integers $\alpha:= i\, i+1$ such that $x_{i+1} = x_i +1$. 
A {\em left ascent} in $x$ is a pair of integers $\alpha := i\, j$ such that 
$1\leq i < j\leq n+1$ and $x_j = x_i+1$.

\begin{Theorem}
Let $W$ denote the symmetric group $S_{n+1}$. 
Then, 
\begin{enumerate}
\item[(1)] the total number of covering relations in $(W,\leq_{LR})$ is $n^2 n!$;
\item[(2)] the number of covering relations of degree 2 in $(W,\leq_{LR})$ is $n n!$.
\end{enumerate}
\end{Theorem}

\begin{proof}
We start with the proof of (2).
Let $x <_{LR} y$ be a covering relation of degree 2 in $S_{n+1}$. 
Then there exist $s_i,s_j\in \{ (1\, 2), (2\, 3),\dots, (n\, n+1)\}$ such that $s_ix = xs_j =y$. 
The left multiplication of $x$ by $s_i$ interchanges the values $x_i$ and $x_{i+1}$ in $x$,
and the right multiplication of $x$ by $s_j$ interchanges the occurrence of $j$ and $j+1$ 
in $x$. Therefore, $x_i = j$ and $x_{i+1}=j+1$. 
Conversely, for each such consecutive pair $x_i x_{i+1}$ in $x=x_1\dots x_{n+1}$ we obtain
a covering relation of degree 2 by interchanging $x_i$ and $x_{i+1}$. 
Therefore, our count is equal to 
\begin{align*}
c_{n+1}:=\text{ the total number of small ascents occurring in permutations in $S_{n+1}$.}
\end{align*}
To find this number let us first fix a small ascent $\alpha = i\, i+1$. 
Clearly, we choose the integer $i$ in $n$ different ways, and $\alpha$ can appear in any of the 
$n!$ permutations of the set $\{1,\dots, i-1, \alpha, i+2,\dots, n+1\}$.
In particular, we see that there are $n \cdot n!$ permutations where $\alpha$ can appear. 
This completes the proof of (2).

Next, we will prove (1). To this end, we will compute 
\begin{align*}
a_{n+1}&:=\text{ the total number of left ascents in $S_{n+1}$,}\\
b_{n+1}&:=\text{ the total number of right ascents in $S_{n+1}$.}
\end{align*}
Then the total number of covering relations is given by $a_{n+1}+b_{n+1}-c_{n+1}$. 
To find $a_{n+1}$, first, choose two positions $i$ and $j$ in $x\in S_{n+1}$,
and set $x_i:= k$ and $x_j:= k+1$ for some $k\in \{1,\dots, n\}$. 
Clearly, there are ${n+1 \choose 2} n$ possible choices. 
Then we choose the remaining entries of $x$ in $(n-1)!$ ways.
Therefore, the total number of left ascents in all permutations in $S_{n+1}$ is given 
by 
$
a_{n+1} = {n+1 \choose 2} n (n-1)! = \frac{n}{2} (n+1)!. 
$
By a similar argument we find that 
$
b_{n+1} =  \frac{n}{2} (n+1)!. 
$
Therefore, we see that 
\begin{align*}
a_{n+1}+b_{n+1}-c_{n+1} = n (n+1)! - n n! = n^2 n!,
\end{align*}
hence, the proof of (1) is complete.
\end{proof}

\textbf{Acknowledgements.} 
This work is partially supported by a grant from the Louisiana Board of Regents.
We thank Kyle Petersen for useful communication about the order complex of the two-sided weak order on $W$.

\bibliography{References}
\bibliographystyle{plain}

\end{document}